\documentclass[%
  a4paper,
  onecolumn,
]{mypreprint}

\usepackage{amsmath}
\usepackage{amssymb}
\usepackage{amsthm}

\usepackage[english]{babel}

\theoremstyle{plain}
\theoremstyle{definition}\newtheorem{remark}{Remark}

\renewcommand{\it}{\itshape}

\usepackage{mathtools}

\usepackage{graphicx}
\usepackage{tikz}
  \usetikzlibrary{arrows}
  \usetikzlibrary{automata}
  \usetikzlibrary{calc}
  \usetikzlibrary{chains}
  \usetikzlibrary{decorations.markings}
  \usetikzlibrary{decorations.pathmorphing}
  \usetikzlibrary{decorations.pathreplacing}
  \usetikzlibrary{fit}
  \usetikzlibrary{matrix}
  \usetikzlibrary{patterns}
  \usetikzlibrary{shapes}

\usepackage{pgfplots}
  \pgfplotsset{compat = 1.13}
  \usepgfplotslibrary{external}
  \tikzset{external/system call = {%
    pdflatex \tikzexternalcheckshellescape
      -halt-on-error
      -interaction=batchmode
      -jobname "\image" "\texsource"}}
  \tikzexternaldisable

\usepackage{caption}
\usepackage{subcaption}

\renewcommand{\rm}[1]{\ensuremath{\mathrm{#1}}}

\newcommand{\trans}{\ensuremath{\mkern-1.5mu\mathsf{T}}}
\newcommand{\herm}{\ensuremath{\mathsf{H}}}
\newcommand{\real}{\ensuremath{\mathrm{Re}}}
\newcommand{\imag}{\ensuremath{\mathrm{Im}}}
\DeclareMathOperator{\dt}{d}
\DeclareMathOperator{\diag}{diag}

\newcommand{\fA}{\ensuremath{\mathcal{A}}}
\newcommand{\fB}{\ensuremath{\mathcal{B}}}
\newcommand{\fC}{\ensuremath{\mathcal{C}}}
\newcommand{\fE}{\ensuremath{\mathcal{E}}}
\newcommand{\fx}{\ensuremath{x}}
\newcommand{\sS}{\ensuremath{\mathcal{S}}}

\newcommand{\fAr}{\ensuremath{\widehat{\fA}}}

\newcommand{\fxr}{\ensuremath{\hat{\fx}}}

\newcommand{\sM}{\ensuremath{M}}
\newcommand{\sD}{\ensuremath{D}}
\newcommand{\sK}{\ensuremath{K}}
\newcommand{\sB}{\ensuremath{B_{\rm{u}}}}
\newcommand{\sC}{\ensuremath{C}}
\newcommand{\sCp}{\ensuremath{\sC_{\rm{p}}}}
\newcommand{\sCv}{\ensuremath{\sC_{\rm{v}}}}
\newcommand{\sx}{\ensuremath{q}}

\newcommand{\sMr}{\ensuremath{\widehat{\sM}}}
\newcommand{\sDr}{\ensuremath{\widehat{\sD}}}
\newcommand{\sKr}{\ensuremath{\widehat{\sK}}}
\newcommand{\sBr}{\ensuremath{\widehat{B}_{\rm{u}}}}
\newcommand{\sCr}{\ensuremath{\widehat{\sC}}}
\newcommand{\sCpr}{\ensuremath{\sCr_{\rm{p}}}}
\newcommand{\sCvr}{\ensuremath{\sCr_{\rm{v}}}}
\newcommand{\sxr}{\ensuremath{\hat{\sx}}}

\renewcommand{\u}{u}
\newcommand{\y}{y}
\newcommand{\yr}{\hat{y}}
\newcommand{\TF}{\ensuremath{H}}
\newcommand{\TFr}{\ensuremath{\widehat{\TF}}}

\newcommand{\ri}{\mathrm{i}}

\newcommand{\R}{\ensuremath{\mathbb{R}}}
\newcommand{\C}{\ensuremath{\mathbb{C}}}
\newcommand{\Hinf}{\ensuremath{\mathcal{H}_{\infty}}}
\newcommand{\Linf}{\ensuremath{\mathcal{L}_{\infty}}}

\newlength{\figurewidthsc}
\newlength{\figurewidthdc}
\newlength{\figureheight}
\usepackage{filemod}

\definecolor{myRed}{HTML}{E34A33}
\definecolor{myBlue}{HTML}{0571B0}
\definecolor{myBrown}{HTML}{A6611A}

\definecolor{plotcolor1}{HTML}{1B9E77}
\definecolor{plotcolor2}{HTML}{D95F02}
\definecolor{plotcolor3}{HTML}{7570B3}
\definecolor{plotcolor4}{HTML}{E7298A}
\definecolor{plotcolor5}{HTML}{66A61E}
\definecolor{plotcolor6}{HTML}{E6AB02}
\definecolor{plotcolor7}{HTML}{A6761D}
\definecolor{plotcolor8}{HTML}{666666}

\definecolor{greedycolor}{rgb}{0.49400,0.18400,0.55600}

\tikzstyle{fom} = [
  line width = 1.2pt,
  black]
\tikzstyle{rom} = [
  line width = 1.2pt,
  myBlue,
  dashed]
\tikzstyle{pole} = [
  line width = 1pt,
  dotted,
  myRed]
\tikzstyle{polefreq} = [
  pole,
  solid,
  mark      = x,
  mark size = 2pt]
\tikzstyle{poleline} = [
  pole,
  solid,
  mark      = o,
  mark size = 2pt]

\tikzstyle{rngline} = [
  black,
  line width = 1.5pt,
  dashed]

\tikzstyle{sobtp} = [
  plotcolor1,
  line width = 1.0pt,
  mark size = 1.5pt,
  mark options = {solid},
  mark = *,
  mark repeat = 20,
  mark phase = 0]

\tikzstyle{sobtv} = [
  plotcolor6,
  line width = 1.0pt,
  mark size = 1.5pt,
  mark options={solid},
  mark=square*,
  mark repeat = 20,
  mark phase  = 5]

\tikzstyle{greedyint} = [
  greedycolor,
  line width = 1.0pt]

\pgfplotscreateplotcyclelist{freqlist}{
  {black, solid, line width = 1.0pt},
  {plotcolor1, line width = 1.0pt, mark size = 1.5pt, mark options = {solid},
    mark = *, mark repeat = 80, mark phase = 0},
  {plotcolor2, line width = 1.0pt, mark size = 1.5pt, mark options={solid},
    mark=o, mark repeat=80, mark phase=8},
  {plotcolor3, line width = 1.0pt, mark size = 1.5pt, mark options={solid},
    mark=x, mark repeat=80, mark phase=16},
  {plotcolor4, line width = 1.0pt, mark size = 1.5pt, mark options={solid},
    mark=+, mark repeat=80, mark phase=24},
  {plotcolor5, line width = 1.0pt, mark size = 1.5pt, mark options={solid},
    mark=triangle*, mark repeat=80, mark phase=32},
  {plotcolor6, line width = 1.0pt, mark size = 1.5pt, mark options={solid},
    mark=square*, mark repeat=80, mark phase=40},
  {plotcolor7, line width = 1.0pt, mark size = 1.5pt, mark options={solid},
    mark=diamond*, mark repeat=80, mark phase=48},
  {plotcolor8, line width = 1.0pt, mark size = 1.5pt, mark options={solid},
    mark=Mercedes star, mark repeat=80, mark phase=56}
}

\pgfplotscreateplotcyclelist{timelist}{
  {black, solid, line width = 1.0pt},
  {plotcolor1, line width = 1.0pt, mark size = 2pt, mark options={solid},
    mark=*, mark repeat=200, mark phase=0},
  {plotcolor2, line width = 1.0pt, mark size = 2pt, mark options={solid},
    mark=o, mark repeat=200, mark phase=20},
  {plotcolor3, line width = 1.0pt, mark size = 2pt, mark options={solid},
    mark=x, mark repeat=200, mark phase=40},
  {plotcolor4, line width = 1.0pt, mark size = 2pt, mark options={solid},
    mark=+, mark repeat=200, mark phase=60},
  {plotcolor5, line width = 1.0pt, mark size = 2pt, mark options={solid},
    mark=triangle*, mark repeat=200, mark phase=80},
  {plotcolor6, line width = 1.0pt, mark size = 2pt, mark options={solid},
    mark=square*, mark repeat=200, mark phase=100},
  {plotcolor7, line width = 1.0pt, mark size = 2pt, mark options={solid},
    mark=diamond*, mark repeat=200, mark phase=120},
  {plotcolor8, line width = 1.0pt, mark size = 2pt, mark options={solid},
    mark=Mercedes star, mark repeat=200, mark phase=140}
}

\begin{document}

\title{Structure-Preserving Model Reduction for Dissipative Mechanical Systems}

\author[$\ast$]{Rebekka S. Beddig}
\affil[$\ast$]{Hamburg University of Technology, Institute of Mathematics,
  Am Schwarzenberg-Campus 3, Geb{\"a}ude~E, 21073 Hamburg, Germany.
  \email{rebekka.beddig@tuhh.de}}

\author[$\ast\ast$]{Peter Benner}
\affil[$\ast\ast$]{Max Planck Institute for Dynamics of Complex Technical
  Systems, Sandtorstra{\ss}e 1, 39106 Magdeburg, Germany.
  \email{benner@mpi-magdeburg.mpg.de} \authorcr \it
  Otto von Guericke University Magdeburg, Faculty of Mathematics,
  Universit{\"a}tsplatz 2, 39106 Magdeburg, Germany.
  \email{peter.benner@ovgu.de}}

\author[$\ast\ast\ast$]{Ines Dorschky}
\affil[$\ast\ast\ast$]{Capgemini Engineering, Olof-Palme-Stra{\ss}e 14,
  81829 M{\"u}nchen, Germany.
  \email{ines.dorschky@uni-hamburg.de}}

\author[$\dagger$]{\authorcr{}Timo Reis}
\affil[$\dagger$]{Technische Universit\"at Ilmenau, Institute of Mathematics,
  Weimarer Stra{\ss}e 25, 98693 Ilmenau, Germany.
  \email{timo.reis@tu-ilmenau.de}}

\author[$\ddag$]{Paul Schwerdtner}
\affil[$\ddag$]{Technische Universit{\"a}t Berlin, Institut f{\"u}r
  Mathematik, Stra{\ss}e des 17. Juni 136, 10623 Berlin, Germany.
  \email{schwerdt@math.tu-berlin.de}}

\author[$\S$]{Matthias Voigt}
\affil[$\S$]{UniDistance Suisse, Schinerstrasse 18, 3900 Brig-VS,
  Switzerland.
  \email{matthias.voigt@fernuni.ch}}

\author[$\P$]{\authorcr{}Steffen W.~R. Werner}
\affil[$\P$]{Courant Institute of Mathematical Sciences, New York University,
  251 Mercer Street, New York, NY 10012, USA.
  \email{steffen.werner@nyu.edu}}

\shorttitle{Structure-Preserving MOR for Dissipative Mechanical Systems}
\shortauthor{R.~S.~Beddig et~al.}
\shortdate{2022-01-26}
\shortinstitute{}


\abstract{%
  Suppressing vibrations in mechanical systems, usually described by
  se\-cond-order dynamical models, is a challenging task in
  mechanical engineering in terms of computational resources even nowadays.
  One remedy is structure-preserving model order reduction to construct
  easy-to-evaluate surrogates for the original dynamical system having the same
  structure.
  In our work, we present an overview of recently developed structure-preserving
  model reduction methods for second-order systems.
  These methods are based on modal and balanced truncation in different
  variants, as well as on rational interpolation.
  Numerical examples are used to illustrate the effectiveness of all described
  methods.
}

\novelty{}

\maketitle


\section{Introduction}%
\label{PBTRMV_sec:intro}

We consider model order reduction of dynamical systems arising from modeling of
mechanical systems, which have the property of \emph{dissipativity}.
That is, energy is only consumed and not produced by the system.
In the particular focus of this work are \emph{linear second-order systems}
\begin{align} \label{PBTRMV_eqn:sosys}
  \begin{aligned}
    \sM \ddot{\sx}(t) + \sD \dot{\sx}(t) + \sK \sx(t) & = \sB \u(t),\\
      \y(t) & = \sCp \sx(t) + \sCv \dot{\sx}(t)
  \end{aligned}
\end{align}
with $\sM, \sD, \sK \in \R^{n \times n}$, $\sB \in \R^{n \times m}$, and
$\sCp, \sCv \in \R^{p \times n}$.
These occur naturally by modeling mechanical systems via force balances, in
which the second derivative of the position vector $\sx(t)$ at time $t\in\R$
occurs by Newton's second law.
Hereby, the matrices $\sM$, $\sD$ and $\sK$ are respectively called
\emph{mass matrix}, \emph{damping matrix} and \emph{stiffness matrix}.
The function $t \mapsto \u(t)$ expresses the \emph{input} to the system (external
forces), a function that can be chosen by the operator (or, alternatively
called, the ``user'') of the system.
Moreover, the model contains an \emph{output} $t \mapsto \y(t)$ that contains
some linear combinations of the state variables and its first derivative, which
are of particular interest.
The typical situation is, especially for systems of high complexity, that the
position vector $\sx(t)$ evolves in a high-dimensional space, that is, the
number $n$ is large.
In contrast to that, the input and output spaces are low-dimensional, i.e.,
$m \ll n$ and $p \ll n$.
Since the number $n$ of position variables is a significant measure for
the difficulty of the numerical simulation of~\cref{PBTRMV_eqn:sosys}, there is
a need for efficient and reliable methods for model reduction, i.e., the
approximation of such systems by ones whose solutions can be computed with
significantly less effort.
In this context, ``reliable'' means that the output of the reduced-order system
is (mathematically proven to be) close to the output of the original system
for the same input signal, whereas ``efficient'' means that the determination of
the reduced-order system comes with as little effort as possible.
Another important demand on model order reduction methods is that they preserve
inherent properties such as stability and the second-order structure of the
system (to mention only a few).
By the latter, we mean that the reduced-order model is of the form
\begin{align} \label{PBTRMV_eqn:sorom}
  \begin{aligned}
    \sMr \ddot{\sxr}(t) + \sDr \dot{\sxr}(t) + \sKr \sxr(t) & = \sBr \u(t),\\
     \yr(t) & = \sCpr \sxr(t) + \sCvr \dot{\sxr}(t)
  \end{aligned}
\end{align}
with $\sMr, \sDr, \sKr \in \R^{r \times r}$, $\sBr \in \R^{r \times m}$
and $\sCpr, \sCvr \in \R^{p \times r}$, and with $r \ll n$.
Moreover, models of mechanical systems have the property that the mass and
stiffness matrices are symmetric positive definite, whereas the negative of the
damping matrix is \emph{dissipative}, that is, $\sD + \sD^{\trans}$ is positive
semi-definite.
These properties are requested to be preserved as well by the reduced-order
system~\cref{PBTRMV_eqn:sorom}.

Meanwhile, model order reduction is an established discipline within applied
mathematics and is subject of textbooks and collections, see~\cite{PBTRMV_Ant05,
PBTRMV_BenMS05, PBTRMV_QuaR14, PBTRMV_BenCOetal17, PBTRMV_BenSGetal21}.
In particular, for first-order systems
\begin{align*}
  \begin{aligned}
    \dot{\fx}(t) & = \fA \fx(t) + \fB \u(t),\\
    \y(t) & = \fC \fx(t),
  \end{aligned}
\end{align*}
there exists a rich theory for their approximation by reduced-order systems of
low state-space dimension; see~\cite{PBTRMV_Ant05} for an overview.
These methods are indeed applicable to first-order representations of
second-order systems like
\begin{align*}
  \begin{aligned}
    \begin{bmatrix} I_{n} & 0 \\ 0 & \sM \end{bmatrix}
      \frac{\dt}{\dt{t}} \begin{bmatrix} \sx(t) \\ \dot{\sx}(t) \end{bmatrix}
      & = \begin{bmatrix} 0 & I_{n} \\ -\sK & -\sD \end{bmatrix}
      \begin{bmatrix} \sx(t) \\ \dot{\sx}(t) \end{bmatrix}
      + \begin{bmatrix} 0 \\ \sB \end{bmatrix} \u(t),\\
    \y(t) & = \begin{bmatrix} \sCp & \sCv \end{bmatrix}
      \begin{bmatrix} \sx(t) \\ \dot{\sx}(t) \end{bmatrix}.
  \end{aligned}
\end{align*}
However, the problem with this is that the reduced-order system is again of
first order, and, in general, it does not have a physical interpretation as a 
mechanical system.
The structure-preserving model order reduction problem of second-order systems
is therefor a problem on its own and new techniques have to be developed.

The model order reduction problem for linear time-invariant systems can also be
considered in the frequency domain.
More precisely, the \emph{transfer function} mapping inputs to putouts in
frequency domain can be considered, which for~\cref{PBTRMV_eqn:sosys} is given
by
\begin{align*}
  \TF(s) & = (\sCp + s \sCv) (s^{2} \sM + s \sD + \sK)^{-1} \sB
    = \begin{bmatrix} \sCp & \sCv \end{bmatrix}
    \begin{bmatrix} sI_{n} & -I_{n} \\ \sK & s \sM + \sD \end{bmatrix}^{-1}
    \begin{bmatrix} 0 \\ \sB \end{bmatrix}.
\end{align*}
Plancherel's theorem~\cite[Prop.~5.1]{PBTRMV_Ant05} provides a link between the
time and frequency domain in a way that -- very roughly speaking -- ``the better
the transfer function of the reduced-order system approximates that of the
original system, the better the outputs of original and reduced-order systems
coincide''.
Two important measures for the distance between transfer functions are the
\emph{$\Hinf$-norm}, which for asymptotically stable systems expresses the
supremal distance between the transfer functions on the imaginary axis; and the
so-called \emph{gap metric}~\cite{PBTRMV_GeoS90}, which applies to arbitrary,
possibly unstable, systems and can be expressed by the $\Hinf$-norm of certain
stable factorizations of transfer functions.
Whereas in the time domain, the $\Hinf$-norm expresses the $L^{2}$-norm
differences of the outputs of the original and reduced-order system, the gap
metric can be seen as a quantitative measure for the distance of the dynamics of
systems.

Besides considering arbitrary linear outputs  $\y(t) = \sCp \sx(t) + \sCv
\dot{\sx}(t)$, in our considerations special emphasis is put on
\emph{co-located} velocity outputs $\y(t) = \sB^{\trans} \dot{\sx}(t)$, which
corresponds to measurements of velocities directly at the force actuators
forming the input.
This special input-output configuration has the additional property that it
provides an energy balance, namely
\begin{align*}
  \forall t \geq 0: \qquad \tfrac12 \left( \dot{\sx}(t)^{\trans} \sM
    \dot{\sx}(t) + \sx(t)^{\trans} \sK \sx(t) \right) - \tfrac12 \left(
    \dot{\sx}(0)^{\trans} \sM \dot{\sx}(0) + \sx(0)^{\trans} \sK \sx(0)
    \right)\\
  = \int\limits_{0}^{t} \y(\tau)^{\trans} \u(\tau) \dt{\tau} -
    \int\limits_{0}^{t} \dot{\sx}(\tau)^{\trans} \sD \dot{\sx}(\tau) \dt{\tau}.
\end{align*}
The expressions $\tfrac12 \dot{\sx}(t)^{\trans} \sM \dot{\sx}(t)$ and $\tfrac12
\sx(t)^{\trans} \sK \sx(t)$ respectively stand for the kinetic and potential
energies of the system at time $t$, whereas
$\int_{0}^{t} \dot{\sx}(\tau)^{\trans} \sD \dot{\sx}(\tau) \dt{\tau}$ is the
dissipated energy, and $\int_{0}^{t} \y(t)^{\trans} \u(t) \dt{t}$ is the energy
put into the system at the actuators within the time interval $[0,t]$.
In particular, since $\sD$ is dissipative and the mass and stiffness matrices
are positive definite, in the case where the system
is in a standstill at $t = 0$, i.e., $\dot{\sx}(0) = \sx(0) = 0$, this energy
balance reduces to
\begin{align*}
  \forall t\geq0:\qquad 0 & \leq \int\limits_{0}^{t} \y(\tau)^{\trans}
    \u(\tau) \dt{\tau}.
\end{align*}
Systems with this property are called \emph{passive}, a property which is
further desired to be preserved by the reduced-order model.
Note that the frequency domain pendant of passivity is \emph{positive realness},
i.e., the transfer function $\TF(s)$ has no poles in the open right complex half
plane and, additionally, $-\TF(s)$ is dissipative in the open right complex half
plane.

For linear time-invariant systems, there are (among others) three ``prominent''
techniques for model order reduction, namely \emph{modal-based},
\emph{balancing-based}, and \emph{interpolation-based} approaches.
The modal-based methods consider eigenvalue problems associated with the
potential poles of the transfer function to retain chosen poles from
the original in the reduced-order model.
Balancing-based methods use energy considerations to figure out parts of the
state only contributing marginally to the input-output behavior, which are
truncated to obtain a reduced-order system of a priori known quality by
providing error bounds, e.g., in the $\Hinf$-norm or gap metric.
The main cost in the determination of reduced-order models by balanced
truncation is the numerical solution of matrix equations of Lyapunov or Riccati
type.
Interpolation-based methods use certain projections of the state space, which
guarantee exactness of the transfer function of the reduced-order system at some
prescribed frequencies.

For second-order systems, the general ideas of balancing-based model order
reduction are subject of various contributions~\cite{PBTRMV_MeyS96,
PBTRMV_ChaLVetal06, PBTRMV_HarVS10, PBTRMV_BenS11, PBTRMV_BenKS13b}; see
also~\cite{PBTRMV_ReiS08} for an overview.
Some progress has been made in preservation of certain physical properties like
passivity in model order reduction of second-order systems with co-located
inputs and outputs \cite{PBTRMV_BenS11, PBTRMV_BenKS13b, PBTRMV_ReiS08}, but
none of these methods are provided with an error bound.
Besides these, there exist interpolatory methods, which succeed either in
preserving the second-order structure~\cite{PBTRMV_BaiS05a, PBTRMV_Fre05,
PBTRMV_SalL06, PBTRMV_BeaG09} or deliver a posteriori $\Hinf$ error
bounds~\cite{PBTRMV_PanWL13}.
However, all the approaches mentioned lack a combination of the two.

In the project ``Structure-Preserving Model Reduction for Dissipative Mechanical
Systems'' of the German Research Foundation (DFG) Priority Program ``Calm,
Smooth, and Smart -- Novel Approaches for Influencing Vibrations by Means of
Deliberately Introduced Dissipation'' (SPP1897), several new approaches for
model order reduction of second-order systems have been developed.
An extract of this work can be found in~\cite{PBTRMV_BedBDetal19}, as well as in
the dissertations of Werner~\cite{PBTRMV_Wer21} and
Dorschky~\cite{PBTRMV_Dor21}.

The structure of this report is as follows.
In \Cref{PBTRMV_sec:somddpa}, we present our results on a dominant pole
algorithm for modally damped mechanical systems.
In \Cref{PBTRMV_sec:so_limited_bt}, a novel balancing-based approach for
second-order systems is presented, which considers the dominant behavior of the
system on some prescribed time and frequency intervals.
In \Cref{PBTRMV_sec:pos_real}, we consider an alternative balancing-based
method for second-order systems with co-located inputs and outputs.
We prove that an error bound in the gap metric holds and that, under some
additional assumptions, a special state-space transformation leads to a
second-order system realization.
\Cref{PBTRMV_sec:hinf_interp} is devoted to an interpolation-based model
order reduction method for second-order systems, based on an optimization-based
technique, which generates a sequence of reduced-order models of descending
error in the $\Hinf$-norm.
The report is concluded in \Cref{PBTRMV_sec:conclusions}.


\section{A dominant pole algorithm for modally damped mechanical systems}%
\label{PBTRMV_sec:somddpa}

One of the oldest model order reduction approaches, which also directly
translates into a structure-preserving setting for second-order
systems~\cref{PBTRMV_eqn:sosys}, is the
method of \emph{modal truncation}~\cite{PBTRMV_Dav66}.
Thereby, the projection basis for the reduced-order model only consists of
the left and right eigenvectors corresponding to the desired eigenvalues.
In case of second-order systems like~\cref{PBTRMV_eqn:sosys}, the corresponding
quadratic eigenvalue problem
\begin{align*}
  \begin{aligned}
    \left( \lambda_{i}^{2} \sM + \lambda_{i} \sD + \sK \right) x_{i} & = 0, &
      y_{i}^{\herm} \left( \lambda_{i}^{2} \sM + \lambda_{i} \sD + \sK \right)
      & = 0
  \end{aligned}
\end{align*}
has to be considered for the left and right eigenvectors
$y_{i}, x_{i} \in \C^{n} \setminus \{0\}$ corresponding to the eigenvalue
$\lambda_{i} \in \C$.
Hereby, $y^{\herm}$ stands for the conjugate transpose of $y$.

For this model order reduction method, the choice of the eigenvalues that
remain in the reduced-order system is critical.
Classical choices are, e.g., taking the rightmost eigenvalues in the complex
plane or the eigenvalues with smallest absolute values.
A significant drawback of those simple choices is the neglection of the input
and output matrices, which have a significant influence on the actual
input-to-output behavior of the system.
The extension of the classical modal truncation method to a more sophisticated
choice of eigenvalues is the
\emph{dominant pole algorithm}~\cite{PBTRMV_MarLP96}.
Here, the eigenvalues with the strongest influence on the system behavior
are computed and then chosen for constructing the reduced-order model.
Adaptations of the dominant pole algorithm to the case of general second-order
systems have been suggested in~\cite{PBTRMV_RomM06} for single-input
single-output systems and in~\cite{PBTRMV_BenKTetal16} for multi-input
multi-output systems.

A modeling approach used very often for mechanical structures results in
modally damped second-order systems.
Hereby, for the second-order system~\cref{PBTRMV_eqn:sosys} it is assumed that
$\sCv = 0$, $\sM, \sD, \sK$ are symmetric positive definite and, additionally,
it holds that $\sD \sM^{-1} \sK = \sK \sM^{-1} \sD$, i.e., the system can be
rewritten into modal coordinates and completely decouples into independent
mechanical systems of order $1$; see, e.g.,~\cite{PBTRMV_BeaB14}.
Classical damping approaches, like Rayleigh and critical damping, fall into
this category.
For this type of mechanical systems, the idea of the dominant pole algorithm
can be reformulated.
As shown in~\cite{PBTRMV_BeaB14}, choosing $X$ as a scaled eigenvector basis
gives
\begin{align} \label{PBTRMV_eqn:eigscale}
  \begin{aligned}
    X^{\trans} \sM X & = \Omega^{-1} && \text{and} &
      X^{\trans} \sK X & = \Omega
  \end{aligned}
\end{align}
with $\Omega = \mathrm{diag}(\omega_{1}, \ldots, \omega_{n}) \in
\R^{n \times n}$ and $X = [x_{1}, \ldots, x_{n}] \in \R^{n \times n}$.
By the modal damping assumption, we further get
\begin{align} \label{PBTRMV_eqn:eigexi}
  X^{\trans} \sD X & = 2\Xi
\end{align}
with $\Xi = \mathrm{diag}(\xi_{1}, \ldots, \xi_{n}) \in \R^{n \times n}$.
Combining~\cref{PBTRMV_eqn:eigscale,PBTRMV_eqn:eigexi}, the
transfer function of~\cref{PBTRMV_eqn:sosys} can be written in a structured
pole-residue form
\begin{align}
  \nonumber
  \TF(s) & = \sCp (s^{2} \sM + s \sD + \sK)^{-1} \sB \\
  \nonumber
    & = \sCp X (s^{2} \Omega^{-1} + 2s \Xi + \Omega)^{-1} X^{\trans} \sB \\
  \label{PBTRMV_eqn:somdtf}
  & = \sum\limits_{k = 1}^{n} \frac{\omega_{k} \left( \sCp x_{k} \right) 
    \left( x_{k}^{\trans} \sB \right)}{(s - \lambda_{k}^{+})(s - 
    \lambda_{k}^{-})},
\end{align}
where the eigenvalue pairs $\lambda_{k}^{+}, \lambda_{k}^{-}$ are given by
\begin{align*}
  \lambda_{k}^{\pm} & = -\omega_{k}\xi_{k} \pm \omega_{k}
    \sqrt{\xi_{k}^{2} - 1}.
\end{align*}
The pole-residue formulation~\cref{PBTRMV_eqn:somdtf} is now used to
derive a new dominant pole algorithm for modally damped second-order systems.
With~\cref{PBTRMV_eqn:somdtf}, an appropriate extension of classical dominant
poles as in~\cite{PBTRMV_BenKTetal16} to pole pairs reads as:
The pole pair $(\lambda_{k}^{+}, \lambda_{k}^{-})$ is called dominant if it
satisfies
\begin{align} \label{PBTRMV_eqn:dompole}
  \frac{ {\lVert \omega_{k}(\sCp x_{k})(x_{k}^{\trans} \sB) \rVert}_{2}}
    {{\lvert \real(\lambda_{k}^{+})(\imag(\lambda_{k}^{+}) \ri -
    \lambda_{k}^{-}) \rvert}}
    & > \frac{ {\lVert\omega_{j}(\sCp x_{j})(x_{j}^{\trans} \sB) \rVert}_{2}}
    {\lvert\real(\lambda_{j}^{+})(\imag(\lambda_{j}^{+}) \ri - \lambda_{j}^{-}) 
    \rvert}
    & \text{for all}~j \neq k.
\end{align}
Note the difference to the dominance measure in~\cite{PBTRMV_SaaSW19}, for
which the two poles of a pair are considered as independent quantities.
The corresponding dominant pole algorithm then computes the $r$ most
dominant pole pairs~\cref{PBTRMV_eqn:dompole} and the corresponding
eigenvectors, such that the reduced-order model's transfer function is given by
\begin{align*}
  \TFr(s) & = \sum\limits_{k = 1}^{r}
    {\frac{\omega_{k} \left( \sCp x_{k} \right) \left( x_{k}^{\trans}
    \sB \right)}
    {(s - \lambda_{k}^{+})(s - \lambda_{k}^{-})}}
    \approx \TF(s),
\end{align*}
with an appropriate ordering of the terms in~\cref{PBTRMV_eqn:somdtf} with
respect to~\cref{PBTRMV_eqn:dompole}.
The projection basis is then given by the $r$ eigenvectors corresponding to
the chosen pole pairs.
The basic ideas of the resulting algorithm are published
in~\cite{PBTRMV_Wer21, PBTRMV_SaaSW19}.
Additionally, we published an implementation of this new algorithm for
large-scale sparse systems as MATLAB and GNU Octave
toolbox~\cite{PBTRMV_BenW21}.

\begin{remark}
  A big advantage of this new approach, compared to the methods
  in~\cite{PBTRMV_RomM06} and~\cite{PBTRMV_BenKTetal16}, is the restriction to
  one-sided projections.
  This preserves the system and eigenvalue structure in each single step such
  that the resulting eigenvector basis will be real and no additional unrelated
  Ritz values, which usually disturb the resulting approximation, are
  introduced in the reduced-order model.
\end{remark}

\begin{figure}[t]
  \centering
  \begin{subfigure}[b]{\textwidth}
    \centering
  \tikzexternalenable%
  \tikzsetnextfilename{butterfly_gyro_tf}%
  \filemodCmp{graphics/butterfly_gyro_tf.tikz}{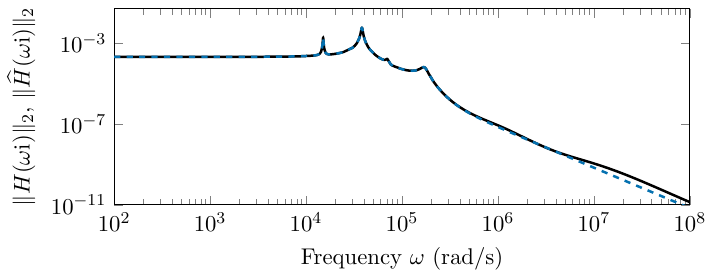}%
    {\tikzset{external/remake next}}{}%
  \begin{tikzpicture}
  \pgfplotstableread{graphics/data/butterfly_gyro_freq.dat}\tableTF
  \pgfplotstableread{graphics/data/butterfly_gyro_pdp.dat}\tablePoles
  
  \begin{loglogaxis}[%
    width  = \figurewidthsc,
    height = \figureheight,
    scale only axis,
    xmin = 1e2,
    xmax = 1e8,
    ymin = 1e-11,
    ymax = 5e-2,
    xminorticks = true,
    yminorticks = false,
    xlabel = {Frequency $\omega$ (rad/s)},
    ylabel = {$\lVert \TF(\omega \ri) \rVert_{2}$, $\lVert \TFr(\omega \ri) \rVert_{2}$}]
    
    \addplot[fom] table[x index=0, y index=1] {\tableTF};
    \addplot[rom] table[x index=0, y index=2] {\tableTF};
  \end{loglogaxis}
\end{tikzpicture}%
  \tikzexternaldisable%

    \caption{Sigma plots.}
  \end{subfigure}

  \begin{subfigure}[b]{\textwidth}
    \centering
  \tikzexternalenable%
  \tikzsetnextfilename{butterfly_gyro_err}%
  \filemodCmp{graphics/butterfly_gyro_err.tikz}{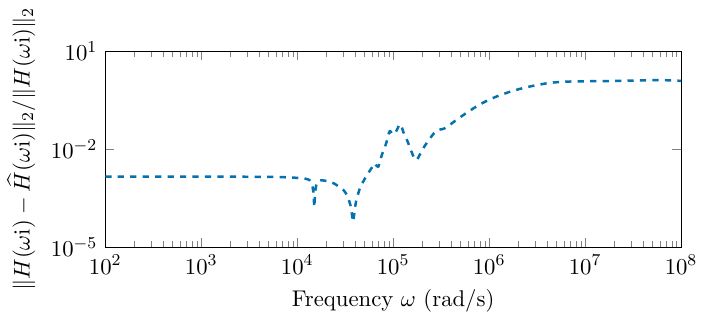}%
    {\tikzset{external/remake next}}{}%
  \begin{tikzpicture}
  \pgfplotstableread{graphics/data/butterfly_gyro_freq.dat}\tableERR
  
  \begin{loglogaxis}[%
    width  = \figurewidthsc,
    height = \figureheight,
    scale only axis,
    xmin = 1e2,
    xmax = 1e8,
    ymin = 1e-5,
    ymax = 1e+1,
    xminorticks = true,
    yminorticks = false,
    xlabel = {Frequency $\omega$ (rad/s)},
    ylabel = {$\lVert \TF(\omega \ri) - \TFr(\omega \ri) \rVert_{2}
      / \lVert \TF(\omega \ri) \rVert_{2}$}]
    
    \addplot[rom] table[x index=0, y index=3] {\tableERR};
  \end{loglogaxis}
\end{tikzpicture}%
  \tikzexternaldisable%

    \caption{Relative error.}
  \end{subfigure}
  \vspace{0\baselineskip}

  \tikzexternalenable%
  \tikzsetnextfilename{butterfly_gyro_legend}%
  \filemodCmp{graphics/butterfly_gyro_legend.tikz}{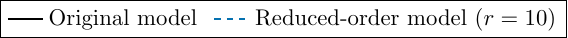}%
    {\tikzset{external/remake next}}{}%
  \begin{tikzpicture}
  \begin{axis}[
    hide axis,
    height = .75\textwidth,
    width  = .12\textheight,
    xmin   = 0,
    xmax   = 10,
    ymin   = 0,
    ymax   = 10,
    legend columns = -1,
    legend style   = {
      anchor = south, 
      /tikz/every even column/.append style={column sep = 2mm}},
        legend cell align = left
  ]
      \addlegendimage{fom}; \addlegendentry{Original model};
      \addlegendimage{rom}; \addlegendentry{Reduced-order model ($r = 10$)};
    \end{axis}
\end{tikzpicture}%
  \tikzexternaldisable%

  \caption{Modally damped dominant pole algorithm results for the butterfly
    gyroscope example.}
  \label{PBTRMV_fig:butterfly_gyro}
\end{figure}

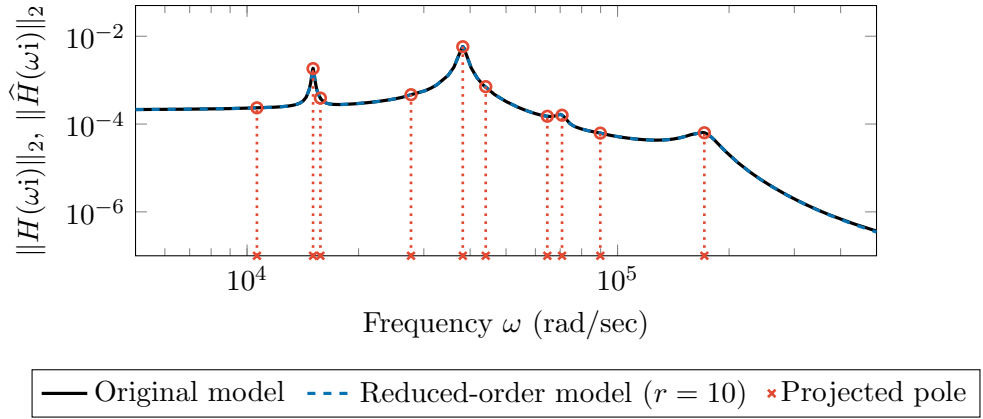
\begin{figure}[t]
  \centering
  \tikzexternalenable%
  \tikzsetnextfilename{butterfly_gyro_pdp}%
  \filemodCmp{graphics/butterfly_gyro_pdp.tikz}{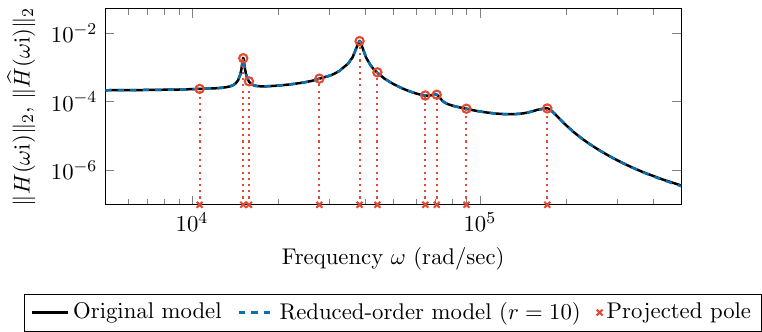}%
    {\tikzset{external/remake next}}{}%
  \begin{tikzpicture}
  \pgfplotstableread{graphics/data/butterfly_gyro_freq_zoom.dat}\tableTF
  \pgfplotstableread{graphics/data/butterfly_gyro_pdp.dat}\tablePoles
  
  \begin{loglogaxis}[%
    width  = \figurewidthsc,
    height = \figureheight,
    scale only axis,
    xmin = 5e3,
    xmax = 5e5,
    ymin = 1e-7,
    ymax = 5e-2,
    xminorticks = true,
    yminorticks = false,
    xlabel = {Frequency $\omega$ (rad/sec)},
    ylabel = {$\lVert \TF(\omega \ri) \rVert_{2}$, $\lVert \TFr(\omega \ri) \rVert_{2}$},
    legend columns = -1,
    legend style   = {
      at = {(.5,-.55)},
      anchor = center,
      /tikz/every even column/.append style = {column sep = 2mm}},
    legend cell align = left
  ]
    
    \addplot[fom] table[x index=0, y index=1] {\tableTF};
    \addlegendentry{Original model};
    
    \addplot[rom] table[x index=0, y index=2] {\tableTF};
    \addlegendentry{Reduced-order model ($r = 10$)};
    
    \addlegendimage{polefreq, only marks};
    \addlegendentry{Projected pole};
    
    \foreach \idx in {0, 2, ..., 18} {
      \addplot[
        pole,
        x filter/.code = {
          \ifnum\coordindex>\numexpr\idx+1\relax
            \def\pgfmathresult{}
          \else
            \ifnum\coordindex<\idx
              \def\pgfmathresult{}
            \fi
          \fi}
      ] table[x index = 0, y index = 1]{\tablePoles};
      
      \addplot[
        polefreq,
        x filter/.code = {
          \ifnum\coordindex>\idx
            \def\pgfmathresult{}
          \else
            \ifnum\coordindex<\idx
              \def\pgfmathresult{}
            \fi
          \fi}
      ] table[x index = 0, y index = 1]{\tablePoles};
      
      \addplot[
        poleline,
        x filter/.code = {
          \ifnum\coordindex>\numexpr\idx+1\relax
            \def\pgfmathresult{}
          \else
            \ifnum\coordindex<\numexpr\idx+1\relax
              \def\pgfmathresult{}
            \fi
          \fi}
      ] table[x index = 0, y index = 1]{\tablePoles};
    }
  \end{loglogaxis}
\end{tikzpicture}%
  \tikzexternaldisable%

  \vspace{0\baselineskip}

  \caption{Projection of the computed dominant poles on the frequency axis
    compared with the transfer functions for the butterfly gyroscope example.}
  \label{PBTRMV_fig:butterfly_gyro_pdp}
\end{figure}

As an illustrative example, we consider the butterfly gyroscope benchmark
from~\cite{PBTRMV_morwiki_gyro} with $n = 17\,361$, $m = 1$ and $p = 12$.
The used Rayleigh damping for the $\sD~(= 10^{-6} \cdot \sK)$ matrix belongs to
the class of modal damping.
We are using the implementation from~\cite{PBTRMV_BenW21} to compute a
reduced-order model with the first $10$ dominant pole pairs by the
criterion~\cref{PBTRMV_eqn:dompole}.
By construction, the resulting reduced-order model has order $10$.
\Cref{PBTRMV_fig:butterfly_gyro} shows the results in the frequency
domain with the frequency response behavior of the original and the
reduced-order model and the pointwise relative error of the approximation.
Up to a frequency of about $10^{6}$\,rad/s, the behavior of the original
system is well reproduced, while later, the two lines begin to diverge slightly.
Additionally, \Cref{PBTRMV_fig:butterfly_gyro_pdp} shows the position of
the computed dominant poles as projection onto the imaginary axis and the
corresponding transfer function behavior.

Comparisons to other model reduction methods and further examples using this
new dominant pole algorithm can be found in~\cite{PBTRMV_Wer21}.

\section{Second-order frequency- and time-limited balanced truncation}%
\label{PBTRMV_sec:so_limited_bt}

A global approximation of the system behavior in frequency or time domain is
often not required in practice.
The second-order limited balanced truncation approaches are a suitable tool
for model order reduction restricted to certain frequency and time ranges.
Thereby, the ideas from the first-order frequency- and time-limited balanced
truncation methods~\cite{PBTRMV_GawJ90} are combined with different second-order
balanced truncation approaches~\cite{PBTRMV_MeyS96, PBTRMV_ChaLVetal06,
PBTRMV_ReiS08}.
A first version of these methods can be found in~\cite{PBTRMV_BedBDetal19,
PBTRMV_SaaSW19, PBTRMV_HaiGIetal18, PBTRMV_HaiGIetal19} and the complete theory
with applications to large-scale sparse systems is contained
in~\cite{PBTRMV_Wer21, PBTRMV_BenW21b}.

The idea of the method is based on the first companion form realization
of~\cref{PBTRMV_eqn:sosys}, which is given by
\begin{align} \label{PBTRMV_eqn:fcsys}
  \begin{aligned}
    \underbrace{\begin{bmatrix} I_{n} & 0 \\ 0 & \sM \end{bmatrix}}_{%
      \phantom{\,\fE}=:\,\fE}
      \dot{\fx}(t) & = \underbrace{\begin{bmatrix} 0 & I_{n} \\ -\sK & -\sD
      \end{bmatrix}}_{\phantom{\,\fA}=:\,\fA} \fx(t) +
      \underbrace{\begin{bmatrix} 0 \\ \sB \end{bmatrix}}_{%
      \phantom{\,\fB}=:\,\fB} \u(t),\\
      \y(t) & = \underbrace{\begin{bmatrix} \sCp & \sCv
      \end{bmatrix}}_{\phantom{\,\fC}=:\,\fC} \fx(t).
  \end{aligned}
\end{align}
For~\cref{PBTRMV_eqn:fcsys}, the classical controllability and observability
Gramians are defined and can be limited as in~\cite{PBTRMV_GawJ90}.
Therefore, the frequency-limited Gramians $P_{\Omega}$ and $Q_{\Omega}$
of~\cref{PBTRMV_eqn:fcsys} are the unique symmetric positive semi-definite
solutions of the potentially indefinite Lyapunov equations
\begin{align*}
  \begin{aligned}
    \fA P_{\Omega} \fE^{\trans} + \fE P_{\Omega} \fA^{\trans} +
      \fB_{\Omega} \fB^{\trans} + \fB \fB_{\Omega}^{\trans} & = 0, \\
    \fA^{\trans} Q_{\Omega} \fE + \fE^{\trans} Q_{\Omega} \fA +
      \fC_{\Omega}^{\trans} \fC + \fC^{\trans}\fC_{\Omega} & = 0,
  \end{aligned}
\end{align*}
for a specified frequency range $\Omega = [\omega_{1}, \omega_{2}] \cup
[-\omega_{2}, -\omega_{1}]$; see~\cite{PBTRMV_GawJ90, PBTRMV_BenKS16}.
The right-hand side matrices contain matrix functions, which are given by
\begin{align*}
  \begin{aligned}
    \fB_{\Omega} & = \real \left( \frac{\ri}{\pi} \ln \left( (\fA + \omega_{2}
      \ri \fE) (\fA + \omega_{1} \ri \fE)^{-1} \right) \right) \fB, \\
    \fC_{\Omega} & = \fC \real \left( \frac{\ri}{\pi} \ln \left( (\fA +
      \omega_{1} \ri \fE)^{-1} (\fA + \omega_{2} \ri \fE) \right) \right).
  \end{aligned}
\end{align*}
Analogously, the time-limited Gramians $P_{\Theta}$ and $Q_{\Theta}$
of~\cref{PBTRMV_eqn:fcsys} are given as the unique symmetric positive
semi-definite solutions of the (potentially) indefinite Lyapunov equations
\begin{align*}
  \begin{aligned}
    \fA P_{\Theta} \fE^{\trans} + \fE P_{\Theta}
      \fA^{\trans} + \fB_{t_{0}}\fB_{t_{0}}^{\trans} -
      \fB_{t_{\rm{f}}} \fB_{t_{\rm{f}}}^{\trans} & = 0,\\
      \fA^{\trans} Q_{\Theta} \fE + \fE^{\trans} Q_{\Theta} \fA +
      \fC_{t_{0}}^{\trans} \fC_{t_{0}} -
      \fC_{t_{\rm{f}}}^{\trans} \fC_{t_{\rm{f}}} & = 0,
  \end{aligned}
\end{align*}
where the time-dependent right-hand sides are defined as
\begin{align*}
  \begin{aligned}
    \fB_{t_{0}} & = \mathrm{e}^{\fA \fE^{-1} t_{0}} \fB, &
      \fB_{t_{\rm{f}}} & = \mathrm{e}^{\fA \fE^{-1}t_{\rm{f}}} \fB, &
      \fC_{t_{0}} & = \fC \mathrm{e}^{\fE^{-1} \fA t_{0}}, &
      \fC_{t_{\rm{f}}} & = \fC \mathrm{e}^{\fE^{-1} \fA t_{\rm{f}}}
  \end{aligned}
\end{align*}
on the time interval $\Theta = [t_{0}, t_{\rm{f}}]$;
see~\cite{PBTRMV_GawJ90, PBTRMV_Kue18}.
Using those Gramians for the different second-order balanced truncation
approaches~\cite{PBTRMV_MeyS96, PBTRMV_ChaLVetal06, PBTRMV_ReiS08} leads to the
second-order limited balanced truncation methods as described
in~\cite{PBTRMV_Wer21, PBTRMV_BenW21b}.
The dense version of the resulting methods is contained in the current version
of the MORLAB toolbox~\cite{PBTRMV_BenW19b, PBTRMV_BenW21c}, and an
implementation for large-scale sparse systems as MATLAB and GNU Octave toolbox
can be found in~\cite{PBTRMV_BenW21a}.

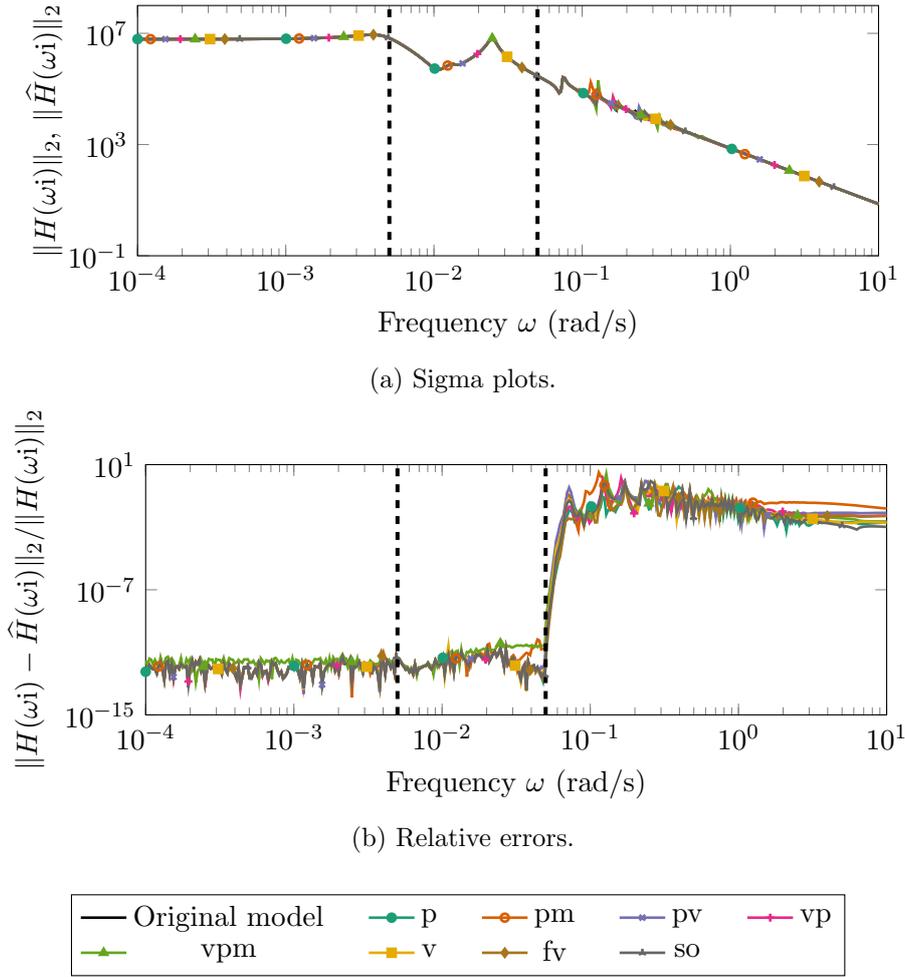
\begin{figure}[t]
  \centering
  \begin{subfigure}[b]{\textwidth}
    \centering
  \tikzexternalenable%
  \tikzsetnextfilename{triplechain_tf}%
  \filemodCmp{graphics/triplechain_tf.tikz}{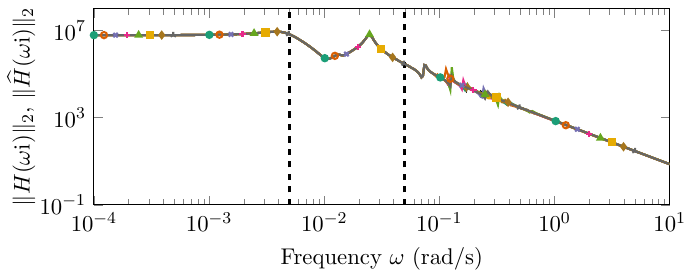}%
    {\tikzset{external/remake next}}{}%
  \begin{tikzpicture}
  \pgfplotstableread{graphics/data/triplechain_freq.dat}\tableTF
  
  \begin{loglogaxis}[%
    width  = \figurewidthsc,
    height = \figureheight,
    scale only axis,
    xmin = 1e-4,
    xmax = 1e+1,
    ymin = 1e-1,
    ymax = 1e+8,
    xminorticks = true,
    xlabel = {Frequency $\omega$ (rad/s)},
    ylabel = {$\lVert \TF(\omega \ri) \rVert_{2}$, $\lVert \TFr(\omega \ri) \rVert_{2}$},
    ylabel style    = { yshift = -.5em},
    draw            = black,
    color           = black,
    cycle list name = freqlist]
    
      \foreach \idx in {1, 2, ..., 9} {
        \addplot table[x index = 0, y index = \idx] {\tableTF};
      }
      
      \addplot[rngline] coordinates {(5e-3, 1e-1) (5e-3, 1e+8)};
      \addplot[rngline] coordinates {(5e-2, 1e-1) (5e-2, 1e+8)};
  \end{loglogaxis}
\end{tikzpicture}
    %
  \tikzexternaldisable%

    \caption{Sigma plots.}
  \end{subfigure}

  \begin{subfigure}[b]{\textwidth}
    \centering
  \tikzexternalenable%
  \tikzsetnextfilename{triplechain_err}%
  \filemodCmp{graphics/triplechain_err.tikz}{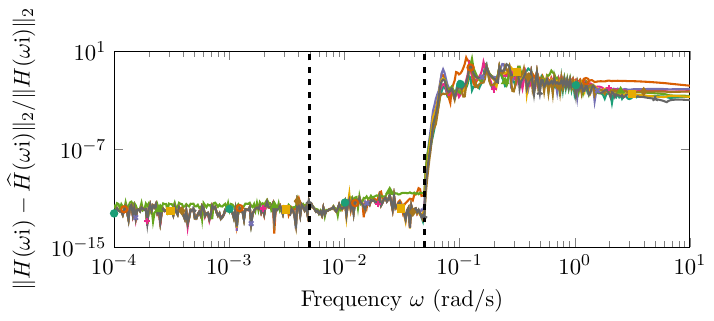}%
    {\tikzset{external/remake next}}{}%
  \begin{tikzpicture}
  \pgfplotstableread{graphics/data/triplechain_freq.dat}\tableERR
  
  \begin{loglogaxis}[%
    width  = \figurewidthsc,
    height = \figureheight,
    scale only axis,
    xmin = 1e-4,
    xmax = 1e+1,
    ymin = 1e-15,
    ymax = 1e+1,
    xminorticks = true,
    xlabel = {Frequency $\omega$ (rad/s)},
    ylabel = {$\lVert \TF(\omega \ri) - \TFr(\omega \ri) \rVert_{2}
      / \lVert \TF(\omega \ri) \rVert_{2}$},
    cycle list name = freqlist]
    
      \pgfplotsset{cycle list shift = 1}
      \foreach \idx in {10, 11, ..., 17} {
        \addplot table[x index = 0, y index = \idx] {\tableERR};
      }
      
      \addplot[rngline] coordinates {(5e-3, 1e-15) (5e-3, 1e+1)};
      \addplot[rngline] coordinates {(5e-2, 1e-15) (5e-2, 1e+1)};
  \end{loglogaxis}
\end{tikzpicture}%
  \tikzexternaldisable%

    \caption{Relative errors.}
  \end{subfigure}
  \vspace{0\baselineskip}

  \tikzexternalenable%
  \tikzsetnextfilename{triplechain_legend}%
  \filemodCmp{graphics/triplechain_legend.tikz}{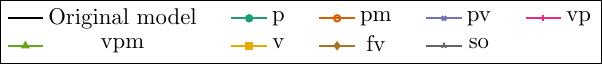}%
    {\tikzset{external/remake next}}{}%
  \begin{tikzpicture}
  \begin{axis}[%
    hide axis,
    scale only axis,
    width = 1mm,
    xmin = -1,
    xmax = 1,
    ymin = -1,
    ymax = 1,
    legend columns  = 5, 
    cycle list name = freqlist,
    legend style = {
      at     = {(0,0)},
      anchor = center,
      /tikz/every even column/.append style = {column sep = 0.5cm}
    }]
      
    \pgfplotsinvokeforeach{1,...,9}{\addplot coordinates {(0,0)};}
        
    \addlegendentry{Original model};
    \addlegendentry{p};
    \addlegendentry{pm};
    \addlegendentry{pv};
    \addlegendentry{vp};
    \addlegendentry{vpm};
    \addlegendentry{v};
    \addlegendentry{fv};
    \addlegendentry{so};
 \end{axis}
\end{tikzpicture}%
  \tikzexternaldisable%

  \caption{Frequency-limited balanced truncation results for the triple chain
    oscillator example.}
  \label{PBTRMV_fig:triplechain}
\end{figure}

As numerical example for the frequency-limited approach, we consider the
triple chain oscillator example as in~\cite{PBTRMV_BedBDetal19}.
We reduce the original model ($n = 1\,201$) by the second-order
frequency-limited balanced truncation method in the interval
$[5 \cdot 10^{-3}, 5 \cdot 10^{-2}]$\,rad/s using the eight different
second-order balancing formulas from~\cite{PBTRMV_BenW21b} to the order
$r = 34$.
The computations are done using the dense implementation of the second-order
frequency-limited balanced truncation method from the latest version
of the MORLAB toolbox~\cite{PBTRMV_BenW19b, PBTRMV_BenW21c}.
The results can be seen in \Cref{PBTRMV_fig:triplechain} with the transfer
functions and the relative approximation errors.
The computed reduced-order models are denoted according to the used balancing
formulas from~\cite{PBTRMV_BenW21b}.
We clearly see the good approximation behavior in the frequency range of
interest.
Only the system computed by the fv formula is stable, while all
others become unstable.
The preservation of stability by the fv formula follows directly from the 
preservation of the symmetric positive definiteness of the system matrices via 
the underlying one-sided projection.
However, note that in general for any balancing formula in the classical
second-order balanced truncation method, there are counter-examples for the
stability preservation~\cite{PBTRMV_ReiS08}.

\begin{figure}[t]
  \centering
  \begin{subfigure}[b]{\textwidth}
    \centering
  \tikzexternalenable%
  \tikzsetnextfilename{singlechain_sim}%
  \filemodCmp{graphics/singlechain_sim.tikz}{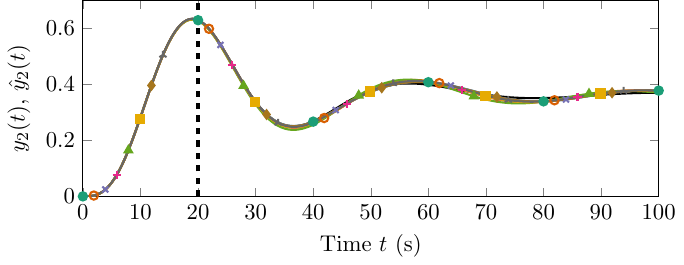}%
    {\tikzset{external/remake next}}{}%
  \begin{tikzpicture}
  \pgfplotstableread{graphics/data/singlechain_time.dat}\tableSIM
  
  \begin{axis}[%
    width  = \figurewidthsc,
    height = \figureheight,
    scale only axis,
    xmin = 0,
    xmax = 100,
    ymin = 0,
    ymax = 0.7,
    xminorticks = true,
    xlabel = {Time $t$ (s)},
    ylabel = {$y_{2}(t)$, $\yr_{2}(t)$},
    draw            = black,
    color           = black,
    cycle list name = timelist]
    
      \foreach \idx in {1, 2, ..., 9} {
        \addplot table[x index = 0, y index = \idx] {\tableSIM};
      }
      
      \addplot[rngline] coordinates {(20, 0) (20, 0.7)};
  \end{axis}
\end{tikzpicture}%
  \tikzexternaldisable%

    \caption{Simulation outputs.}
  \end{subfigure}

  \begin{subfigure}[b]{\textwidth}
    \centering
  \tikzexternalenable%
  \tikzsetnextfilename{singlechain_err}%
  \filemodCmp{graphics/singlechain_err.tikz}{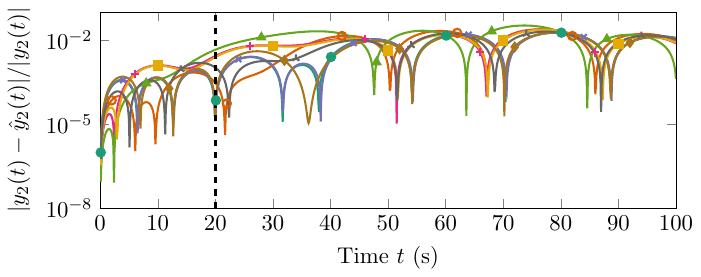}%
    {\tikzset{external/remake next}}{}%
  \begin{tikzpicture}
  \pgfplotstableread{graphics/data/singlechain_time.dat}\tableERR
  
  \begin{semilogyaxis}[%
    width  = \figurewidthsc,
    height = \figureheight,
    scale only axis,
    xmin = 0,
    xmax = 100,
    ymin = 1e-8,
    ymax = 1e-1,
    xminorticks = true,
    xlabel = {Time $t$ (s)},
    ylabel = {$\lvert y_{2}(t) - \yr_{2}(t) \rvert / \lvert y_{2}(t) \rvert$},
    draw            = black,
    color           = black,
    cycle list name = timelist]
    
    \pgfplotsset{cycle list shift = 1}
    \foreach \idx in {10, 11, ..., 17} {
      \addplot table[x index = 0, y index = \idx] {\tableERR};
    }
      
    \addplot[rngline] coordinates {(20, 1e-8) (20, 1e-1)};
  \end{semilogyaxis}
\end{tikzpicture}%
  \tikzexternaldisable%

    \caption{Relative errors.}
  \end{subfigure}
  \vspace{0\baselineskip}

  \tikzexternalenable%
  \tikzsetnextfilename{singlechain_legend}%
  \filemodCmp{graphics/singlechain_legend.tikz}{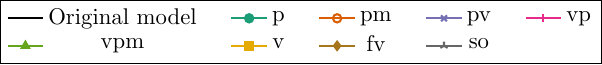}%
    {\tikzset{external/remake next}}{}%
  \begin{tikzpicture}
  \begin{axis}[%
    hide axis,
    scale only axis,
    width = 1mm,
    xmin = -1,
    xmax = 1,
    ymin = -1,
    ymax = 1,
    legend columns  = 5, 
    cycle list name = timelist,
    legend style = {
      at     = {(0,0)},
      anchor = center,
      /tikz/every even column/.append style = {column sep = 0.5cm}
    }]
      
    \pgfplotsinvokeforeach{1,...,9}{\addplot coordinates {(0,0)};}
        
    \addlegendentry{Original model};
    \addlegendentry{p};
    \addlegendentry{pm};
    \addlegendentry{pv};
    \addlegendentry{vp};
    \addlegendentry{vpm};
    \addlegendentry{v};
    \addlegendentry{fv};
    \addlegendentry{so};
 \end{axis}
\end{tikzpicture}%
  \tikzexternaldisable%

  \caption{Time-limited balanced truncation results for the single chain
    oscillator example.}
  \label{PBTRMV_fig:singlechain}
\end{figure}

To illustrate the time-limited balanced truncation method, we use the single
chain oscillator example as described in~\cite{PBTRMV_BenW21b}.
Here, we use the implementation for large-scale sparse mechanical systems
from~\cite{PBTRMV_BenW21a} to reduce the original system ($n = 12\,000$) to
order $r = 3$ in the time interval $[0, 20]$\,s.
The results for the second output entry can be seen in
\Cref{PBTRMV_fig:singlechain}, where in the time region of interest, the
original system is nicely approximated by the reduced-order models.
For all balancing formulas, the resulting systems are stable.

More detailed comparisons of the different balancing formulas and further
examples can be found in~\cite{PBTRMV_Wer21, PBTRMV_BenW21b}.


\section{Positive real balanced truncation for second-order systems}%
\label{PBTRMV_sec:pos_real}

In this section, we consider second-order systems of the
form~\cref{PBTRMV_eqn:sosys}, where $\sM, \sK > 0$, $\sD \geq 0$, and either
we exclusively measure positions, i.e., $\sCv = 0$, or velocities, i.e.,
$\sCp = 0$.
We start with the second case and additionally assume co-located inputs and
outputs, which means $\sB = \sCv^{\trans}$.
This case is treated in~\cite{PBTRMV_DorRV21}.
Using the Cholesky factorizations $\sK = G G^{\trans}$ and $M = N N ^{\trans}$,
the system can be rewritten in first-order form as
\begin{align} \label{PBTRMV_eqn:pfosys}
  \begin{aligned}
    \dot{\fx}(t) & = \underbrace{\begin{bmatrix} 0 & G^{\trans} N^{-\trans} \\
      -N^{-1} G & -N^{-1} \sD N^{-\trans}
      \end{bmatrix}}_{\phantom{\,\fA}=:\,\fA} \fx(t) +
      \underbrace{\begin{bmatrix} 0 \\ N^{-1} \sB
      \end{bmatrix}}_{\phantom{\,\fB}=:\,\fB} \u(t),\\
    \y(t) & = \underbrace{\begin{bmatrix} 0 & \sB^{\trans} N^{-\trans}
      \end{bmatrix}}_{\phantom{\,\fC}=:\,\fC} \fx(t).
  \end{aligned}
\end{align}
Besides passivity, the most important feature of this system is that it has an
internal symmetry structure $\fA\sS_{n}=\sS_{n} \fA^{\trans}$ and $\fC = 
\fB^{\trans} = \fB^{\trans} \sS_n$, where $\sS_{n}:=\diag(-I_{n},I_{n})$.
In particular, its transfer function $\TF(s) =\fC(sI_{2n}-\fA)^{-1}\fB$ is
symmetric, i.e., it fulfills $\TF(s)^{\trans} = \TF(s)$.
We will make heavy use of this symmetry structure.

The model reduction technique consists of two steps:

\textit{Step~1:} We apply positive real balanced
truncation~\cite{PBTRMV_GuiO13} to the first-order
system~\cref{PBTRMV_eqn:pfosys}.
Here, the computational bottleneck is determining the numerical solution of
the Lur'e equations
\begin{align*}
  \begin{aligned}
    \fA^{\trans} P + P \fA & = -K_{\rm{c}}^{\trans} K_{\rm{c}}, &
      \fA Q + Q \fA^{\trans} & = -K_{\rm{o}}^{\trans} K_{\rm{o}},\\
    P \fB - \fC^{\trans} & = 0, &
      Q \fC^{\trans} - \fB & = 0,
  \end{aligned}
\end{align*}
for stabilizing solutions $P, Q \in \R^{2n \times 2n}$; see \cite{PBTRMV_Rei11}.
The internal symmetry structure of \cref{PBTRMV_eqn:pfosys} yields
$Q = \sS_{n} P \sS_{n}$, hence only the Lur'e equation $\fA^{\trans} P + P\fA =
-K_{\rm{c}}^{\trans} K_{\rm{c}}$, $P \fB - \fC^{\trans} = 0$ has to be solved
for the matrices $P$ and $K_{\rm{c}}$.
We can use the method from \cite{PBTRMV_PolR12} to obtain a low-rank
approximate solution $P \approx L^{\trans} L$.
The sign symmetry of the first-order system~\cref{PBTRMV_eqn:pfosys} yields
that its positive real characteristic values (which are defined to be the
positive square roots of the eigenvalues of $P Q$) can in a certain sense be
allocated to the symmetry structure of the system~\cref{PBTRMV_eqn:pfosys}.
More precisely, the positive real characteristic values are the absolute values
of the eigenvalues of the symmetric matrix $L\sS_{n} L^{\trans}$.
By truncating equally many states corresponding to positive and negative
eigenvalues of $L\sS_{n} L^{\trans}$, it is shown that the resulting
first-order model is -- without any further computational effort -- of the form
\begin{align} \label{PBTRMV_eqn:balanced_so0}
  \begin{aligned}
    \dot{\fxr}(t) & = 
      \begin{bmatrix}
      0 & 0 & 0 & 0 & 0& \fA_{16} \\
      0  & 0 & 0 &0 & \fA_{25} & \fA_{26} \\
      0 & 0 & \fA_{33} & \fA_{34} & 0& \fA_{36} \\
      0 & 0 & -\fA_{34}^{\trans} & \fA_{44} & 0 & \fA_{46} \\
      0 & -\fA_{25}^{\trans} & 0 & 0 & 0 & 0 \\
      -\fA_{16}^{\trans} & -\fA_{26}^{\trans} & -\fA_{36}^{\trans} &
      \fA_{46}^{\trans} & 0 & \fA_{66}
      \end{bmatrix}{\fxr}(t)+\begin{bmatrix}
      0 \\ 0 \\ 0 \\ 0 \\ 0 \\ \fB_6
      \end{bmatrix}u(t),\\
    \yr(t) & = \begin{bmatrix}
      0 & 0 & 0 & 0 & 0 & \fB_6^{\trans}
      \end{bmatrix}\fxr(t),
  \end{aligned}
\end{align}
where the block sizes from left to right and from top to bottom are
$m, \ell, p, p, \ell, m$, with $r = p + m + \ell$.
Note that, if $\fA_{33}$ is zero, then -- by merging some of the blocks --
\cref{PBTRMV_eqn:balanced_so0}~has the form
\begin{align} \label{PBTRMV_eqn:redsys1}
  \begin{aligned}
    \dot{\fxr}(t)
      & = \begin{bmatrix}
      0 & \widehat{G}^{\trans}  \\ -\widehat{G} & -\widehat{D}
      \end{bmatrix} \fxr(t)+
      \begin{bmatrix} 0 \\ \sBr \end{bmatrix} u(t),\\
    \yr(t) & = \begin{bmatrix} 0 & \sBr^{\trans} \end{bmatrix}\fxr(t),
  \end{aligned}
\end{align}
which results in a reduced-order model in second-order
form~\cref{PBTRMV_eqn:sorom} with $\sMr = I_{r}$,
$\sKr = \widehat{G} \widehat{G}^{\trans}$, 
$\sCpr = 0$, and $\sCvr = \sBr^{\trans}$.
This is regrettably not the case in general, which is why we apply
the following step.

\textit{Step~2:} We apply a state-space transformation
to~\cref{PBTRMV_eqn:balanced_so0} such that the matrix $\fA_{33}$ vanishes.
More precisely, we first intend to find some invertible $T\in\R^{2p\times2p}$
that preserves the symmetry structure, i.e., it fulfills
$T^{\trans} \sS_{p} T = \sS_{p}$ and
\begin{align*}
  T^{-1} \begin{bmatrix} \fA_{33} & \fA_{34} \\ -\fA_{34}^{\trans} & \fA_{44}
  \end{bmatrix} T & =\begin{bmatrix} 0 & \fAr_{34} \\
  -\fAr_{34}^{\trans} & \fAr_{44} \end{bmatrix}.
\end{align*}
Then, a state-space transformation with $V = \diag(I_{m+\ell}, T, I_{\ell+m})$
leads to a system which is indeed of the form~\cref{PBTRMV_eqn:redsys1} and
can then be rewritten as a second-order system.
Such a transformation is based on techniques from indefinite linear
algebra~\cite{PBTRMV_GohLR83}, and can be computed without any remarkable
computational effort.
The existence of such a transformation is linked to the $2k$ real eigenvalues
of $\begin{bsmallmatrix} \fA_{33} & \fA_{34} \\ -\fA_{34}^{\trans} & \fA_{44}
\end{bsmallmatrix}$.
Under the assumption of semi-simplicity, these eigenvalues can be assigned a
certain signature structure according to the sign of
$-v_{i}^{\trans} \sS_{k} v_{i}$, where $v_{i}$ is the eigenvector of the
$i$-th eigenvalue.
In particular, we have $k$ eigenvalues $\mu_{1}^{-} \le \ldots \le \mu_{k}^{-}$
of negative type and $k$ eigenvalues $\mu_{1}^{+} \le \ldots \le\mu_{k}^{+}$ of
positive type. 
It is then shown that such a transformation is possible, if
$\mu_{i}^{-} < \mu_{i}^{+}$ for all $i = 1, \ldots, k$.
A sufficient criterion on the original system for the existence of such a
transformation is that it is \emph{overdamped}, that is $\big(v^{\trans} \sD
v \big)^2 > 4 \big(v^{\trans} \sM v\big) \big(v^{\trans} \sK v \big)$
for all $v \in \R^{n}$.

The resulting second-order system in particular fulfills $\sMr, \sKr > 0$,
and it is further shown that $\sDr = \sDr^{\trans}$ has at most $m$ negative
eigenvalues.
If the original system is overdamped, then even $\sDr > 0$. 
Moreover, the gap metric distance between the transfer functions $\TF(s)$ of
the original and $\TFr(s)$ of the reduced-order system is shown to be bounded
from above by twice the sum of the truncated positive real characteristic
values. 

Considering a system dilation, we are able to extend the previously presented
reduction to second-order systems with velocity measurements $\y(t) = 
\sCp \sx(t)$, which are not necessarily co-located to the input.
This can be done by considering the extended system
\begin{align*}
  \begin{aligned}
    \sM \ddot{\sx}(t) + \sD \dot{\sx}(t) + \sK \sx(t) & =
      \begin{bmatrix} \sB & \sCv^{\trans} \end{bmatrix}
      \begin{bmatrix} \u_{1}(t) \\ \u_{2}(t) \end{bmatrix}, \\
    \begin{bmatrix} \y_{1}(t) \\ \y_{2}(t) \end{bmatrix} & =
      \begin{bmatrix} \sB^{\trans} \\ \sCv \end{bmatrix} \dot{\sx}(t).
  \end{aligned}
\end{align*}
This system is again positive real and from the algorithm above we thus obtain a
reduced-order system for the extended system as
\begin{align} \label{PBTRMV_eqn:redsysnoncolredext}
  \begin{aligned}
    \sMr \ddot{\sxr}(t) + \sDr \dot{\sxr}(t) + \sKr \sxr(t) & =
      \begin{bmatrix} \sBr & \sCvr^{\trans} \end{bmatrix}
      \begin{bmatrix} \u_{1}(t) \\ \u_{2}(t) \end{bmatrix},\\
    \begin{bmatrix} \yr_1(t) \\ \yr_2(t) \end{bmatrix} & =
      \begin{bmatrix} \sBr^{\trans} \\ \sCvr
    \end{bmatrix} \dot{\sxr}(t),
  \end{aligned}
\end{align}
where $\sBr \in \R^{r \times m}$ and $\sCvr \in \R^{p \times r}$.
From that we obtain a reduced-order system for~\cref{PBTRMV_eqn:sosys} as
\begin{align} \label{PBTRMV_eqn:redsysnoncolred}
  \begin{aligned}
    \sMr \ddot{\sxr}(t) + \sDr \dot{\sxr}(t) + \sKr \sxr(t) & = \sBr
      \u_{1}(t),\\
    \yr_{2}(t) & = \sCvr^{\trans} \dot{\sxr}(t).
  \end{aligned}
\end{align}
As the transfer function of~\cref{PBTRMV_eqn:redsysnoncolred} is a submatrix
of the transfer function of the system~\cref{PBTRMV_eqn:redsysnoncolredext},
and the same holds for the original system, the reduction error in the gap
metric is bounded by the one of the extended system.

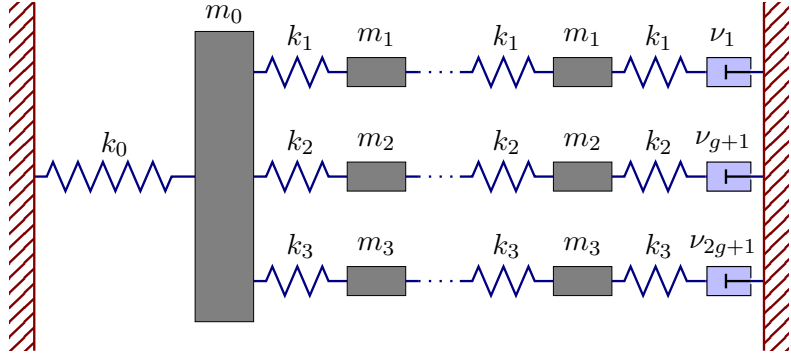
\begin{figure}[t]
  \centering
  \tikzexternalenable%
  \tikzsetnextfilename{mass_spring_damper_3n_TruV09}%
  \filemodCmp{graphics/mass_spring_damper_3n_TruV09.tikz}{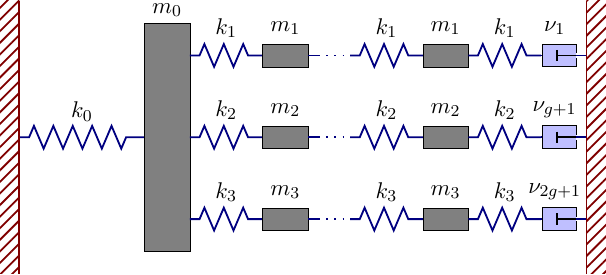}%
    {\tikzset{external/remake next}}{}%
%
%
%

\begin{tikzpicture}[
  x = 2.845em,
  y = 2.845em,
  inner sep = 0,
  outer sep = 0,
  draw      = black]
  
  \colorlet{springcolor}{blue!50!black}
  \colorlet{groundcolor}{red!50!black}
  \colorlet{dampercolor}{blue!25}
  \colorlet{masscolor}{black!50}

  \tikzstyle{every node} = [color = black]
  
  \makeatletter
  \tikzset{
    hatch distance/.store in = \hatchdistance,
    hatch distance = .3em,
    hatch thickness/.store in = \hatchthickness,
    hatch thickness = .08em}
  \pgfdeclarepatternformonly[\hatchdistance,\hatchthickness]{north east hatch}
    {\pgfpoint{-.1em}{-.1em}}
    {\pgfpoint{\hatchdistance+.5*\hatchthickness}%
      {\hatchdistance+.5*\hatchthickness}}
    {\pgfpoint{\hatchdistance}{\hatchdistance}}{
      \pgfsetcolor{\tikz@pattern@color}
      \pgfsetlinewidth{\hatchthickness}
      \pgfpathmoveto{\pgfpoint{0em}{0em}}
      \pgfpathlineto{\pgfpoint{\hatchdistance}{\hatchdistance}}
      \pgfusepath{stroke}}
  \makeatother
  
  \tikzstyle{damper} = [
    line width = .08em,
    draw       = springcolor,
    decoration = {
      markings,  
      mark connection node = dmp,
      mark = at position 0.4 with {
        \node (dmp) [
          inner sep      = 0em,
          minimum height = 1em,
          minimum width  = 1.5em,
          draw = none,
          fill = dampercolor] {};
        \draw [
          line width = .04em,
          draw       = black
        ] ($(dmp.east) + (0, .1em)$)
          -- (dmp.north east)
          -- (dmp.north west)
          -- (dmp.south west)
          -- (dmp.south east)
          -- ($(dmp.east) - (0, .1em)$);
        \draw [
          line width = .08em,
          draw = black
        ] (dmp.east) -- ($(dmp.center) + (-.1em,0)$) coordinate (dmp2);
        \draw [
          line width = .08em,
          draw = black
        ] ($(dmp2) + (0,.25em)$) -- ($(dmp2) - (0,.25em)$);
      }}, 
    decorate
  ]
  
  \tikzstyle{ground} = [
    fill,
    pattern         = north east hatch,
    hatch distance  = .51em,
    hatch thickness = .08em,
    pattern color   = groundcolor,
    draw            = none,
    minimum width   = 12em,
    minimum height  = .8535em,
    rotate          = -90
  ]
  
  \tikzstyle{mass} = [
    rectangle,
    minimum height = 1em,
    minimum width  = 2em,
    line width     = .02em,
    fill           = masscolor,
    draw
  ]
  
  \tikzstyle{spring} = [
    line width = .08em,
    springcolor,
    decorate, 
    decoration = {
      zigzag, 
      pre length     = .4268em, 
      post length    = .4268em, 
      segment length = .8536em,
      amplitude      = .5em}
  ]
  
  \node[ground] (wall1) {};

  \draw[spring] (wall1.north) -- ++(5.5em,0)
    node[midway, above = .7em, draw = none] {$k_{0}$}
    coordinate (tmp);
  \node[mass, minimum height = 10em, label = {[label distance = .3em]:$m_{0}$}]
    (mass0) at (tmp) [anchor = west] {};
  
  \coordinate (txtline) at ($(mass0.east) + (0,4.8em)$);
  \draw[spring] ($(mass0.east) + (0,3.6em)$) -- ++(3.2em,0)
    coordinate[midway] (txt)
    coordinate (tmp);
  \node[draw = none] at (txt |- txtline) {$k_{1}$};
  \node[mass, anchor = west] (mass) at (tmp) {};
  \node[draw = none] at (mass.center |- txtline) {$m_{1}$};
  \draw[line width = .08em, springcolor]
    (mass.east) -- ++(.4268em,0) coordinate (tmp);
  \draw[line width = .08em, springcolor,
    dash pattern = on \pgflinewidth off .2845em]
    (tmp) -- ++(1.4225em,0) coordinate (tmp);
  \draw[spring] (tmp) -- ++(3.2em,0)
    coordinate[midway] (txt)
    coordinate (tmp);
  \node[draw = none] at (txt |- txtline) {$k_{1}$};
  \node[mass, anchor = west] (mass) at (tmp) {};
  \node[draw = none] at (mass.center |- txtline) {$m_{1}$};
  \draw[spring] (mass.east) -- ++(3.2em,0)
    coordinate[midway] (txt)
    coordinate (tmp);
  \node[draw = none] at (txt |- txtline) {$k_{1}$};
  \draw[damper] (tmp) -- ++(2em,0)
    coordinate[pos = .3] (txt)
    coordinate (tmp);
  \node[draw = none] at (txt |- txtline) {$\nu_{1}$};
  
  \coordinate (txtline) at ($(mass0.east) + (0,1.2em)$);
  \draw[spring] (mass0.east) -- ++(3.2em,0)
    coordinate[midway] (txt)
    coordinate (tmp);
  \node[draw = none] at (txt |- txtline) {$k_{2}$};
  \node[mass, anchor = west] (mass) at (tmp) {};
  \node[draw = none] at (mass.center |- txtline) {$m_{2}$};
  \draw[line width = .08em, springcolor]
    (mass.east) -- ++(.4268em,0) coordinate (tmp);
  \draw[line width = .08em, springcolor,
    dash pattern = on \pgflinewidth off .2845em]
    (tmp) -- ++(1.4225em,0) coordinate (tmp);
  \draw[spring] (tmp) -- ++(3.2em,0)
    coordinate[midway] (txt)
    coordinate (tmp);
  \node[draw = none] at (txt |- txtline) {$k_{2}$};
  \node[mass, anchor = west] (mass) at (tmp) {};
  \node[draw = none] at (mass.center |- txtline) {$m_{2}$};
  \draw[spring] (mass.east) -- ++(3.2em,0)
    coordinate[midway] (txt)
    coordinate (tmp);
  \node[draw = none] at (txt |- txtline) {$k_{2}$};
  \draw[damper] (tmp) -- ++(2em,0)
    coordinate[pos = .3] (txt)
    coordinate (tmp);
  \node[draw = none] at (txt |- txtline) {$\nu_{g+1}$};
  
  \coordinate (txtline) at ($(mass0.east) - (0,2.4em)$);
  \draw[spring] ($(mass0.east) - (0,3.6em)$) -- ++(3.2em,0)
    coordinate[midway] (txt)
    coordinate (tmp);
  \node[draw = none] at (txt |- txtline) {$k_{3}$};
  \node[mass, anchor = west] (mass) at (tmp) {};
  \node[draw = none] at (mass.center |- txtline) {$m_{3}$};
  \draw[line width = .08em, springcolor]
    (mass.east) -- ++(.4268em,0) coordinate (tmp);
  \draw[line width = .08em, springcolor,
    dash pattern = on \pgflinewidth off .2845em]
    (tmp) -- ++(1.4225em,0) coordinate (tmp);
  \draw[spring] (tmp) -- ++(3.2em,0)
    coordinate[midway] (txt)
    coordinate (tmp);
  \node[draw = none] at (txt |- txtline) {$k_{3}$};
  \node[mass, anchor = west] (mass) at (tmp) {};
  \node[draw = none] at (mass.center |- txtline) {$m_{3}$};
  \draw[spring] (mass.east) -- ++(3.2em,0)
    coordinate[midway] (txt)
    coordinate (tmp);
  \node[draw = none] at (txt |- txtline) {$k_{3}$};
  \draw[damper] (tmp) -- ++(2em,0)
    coordinate[pos = .3] (txt)
    coordinate (tmp);
  \node[draw = none] at (txt |- txtline) {$\nu_{2g+1}$};
  
  \node[ground, anchor = south] (wall2) at (tmp |- wall1.north) {};
  
  \draw[line width = .08em, groundcolor]
    (wall1.north west) -- (wall1.north east);
  \draw[line width = .08em, groundcolor]
    (wall2.south west) -- (wall2.south east);
\end{tikzpicture}%
  \tikzexternaldisable%

  \caption{(3g + 1)-mass (triple chain) oscillator with three
    dampers~\cite{PBTRMV_TruV09}.}
  \label{PBTRMV_fig:triplechain1}
\end{figure}

We illustrate the performance of the reduction method above with an example of
three coupled mass-spring-damper chains; see~\cite[Ex.~2]{PBTRMV_TruV09}.
The triple chain consists of three rows that are coupled via a mass $m_{0}$,
which is connected to the fixed base with a spring with stiffness $k_{0}$.
Each row contains $g$ masses, $g+1$ springs and one damper, which is attached to
a wall; see \Cref{PBTRMV_fig:triplechain1}.
One can write the free system as
\begin{align*}
  \sM \ddot{\sx}(t) + \sD \dot{\sx}(t) + \sK \sx(t) & = 0,
\end{align*}
where $\sM, \sD$, and $\sK$ are defined as $\sM = \diag(m_{1},\ldots, m_{1},
m_{2}, \ldots, m_{2},\, m_{3}, \ldots,  m_{3})$, $\sD = \alpha \sM +
\beta \sK + \nu_{1} e_{1} e_{1}^{\trans} + \nu_{g+1} e_{g+1} e_{g+1}^{\trans} +
\nu_{2g+1} e_{2g+1} e_{2g+1}^{\trans}$ and
\begin{align*}
  \begin{aligned}
    \sK & = \begin{bmatrix}
      K_{11} & & & -\kappa_{1} \\
      & K_{22} & & -\kappa_{2} \\
      & & K_{33} & -\kappa_{3} \\
      -\kappa_{1}^{\trans} & -\kappa_{2}^{\trans} & -\kappa_{3}^{\trans} &
      ~ k_{1} + k_{2} + k_{3} + k_{0}
      \end{bmatrix}, &
      K_{ii} & = k_i\begin{bmatrix}
      2 & -1 \\ -1 & 2 & -1 \\
      & \ddots & \ddots & \ddots \\
      & & -1 & 2 & -1 \\
      & & & -1 & 2
      \end{bmatrix}
  \end{aligned}
\end{align*}
with $\kappa_{i} = \begin{bmatrix} 0 & \ldots & 0 & k_{i} \end{bmatrix}^{\trans}
\in \R^{1\times g}$ and $K_{ii} \in \R^{g \times g}$ for $i = 1, 2, 3$.
We choose the input $\sB = \begin{bmatrix} 1 & \ldots & 1
\end{bmatrix}^{\trans}$ and equally measure the velocities such that
$\sCv = \sB^{\trans}$.
The second-order control system reads
\begin{align*}
  \begin{aligned}
    \sM \ddot{\sx}(t) + \sD \dot{\sx}(t) + \sK \sx(t) & = \sB \u(t), & 
    \y(t) & = \sB^{\trans} \dot{\sx}(t).
  \end{aligned}
\end{align*}
We consider the triple chain with $g = 500$, thus the number of differential
equations is $n = 3g + 1 = 1\,501$.
The stiffness and mass parameters are set as
\begin{align*}
  \begin{aligned}
    k_{0} & = 50, & k_{1} & = 10, & k_{2} & = 20, & k_{3} & = 1,\\
    m_{0} & = 1, & m_{1} & = 1, & m_{2} & = 2, & m_{3} & = 3,
  \end{aligned}
\end{align*}
and the damping parameters $\alpha = \beta = 0.002$ and $\nu_{1} =
\nu_{g+1} = \nu_{2g+1} = 5$.

\begin{figure}[t]
  \centering
  \begin{subfigure}{\textwidth}
    \centering
  \tikzexternalenable%
  \tikzsetnextfilename{soprbt_systems}%
  \filemodCmp{graphics/soprbt_systems.tikz}{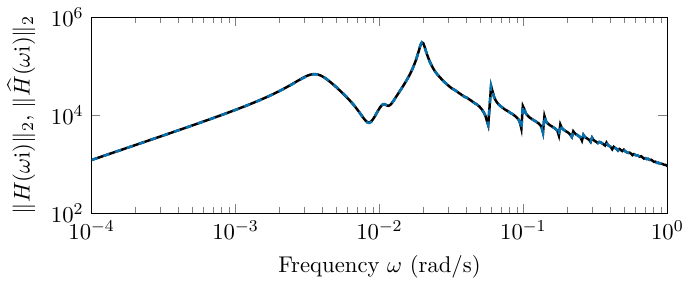}%
    {\tikzset{external/remake next}}{}%
  \input{graphics/soprbt_systems.tikz}%
  \tikzexternaldisable%

    \caption{Sigma plots.}
    \label{PBTRMV_fig:soprbt_tf}
  \end{subfigure}

  \begin{subfigure}{\textwidth}
    \centering
  \tikzexternalenable%
  \tikzsetnextfilename{soprbt_errorrel}%
  \filemodCmp{graphics/soprbt_errorrel.tikz}{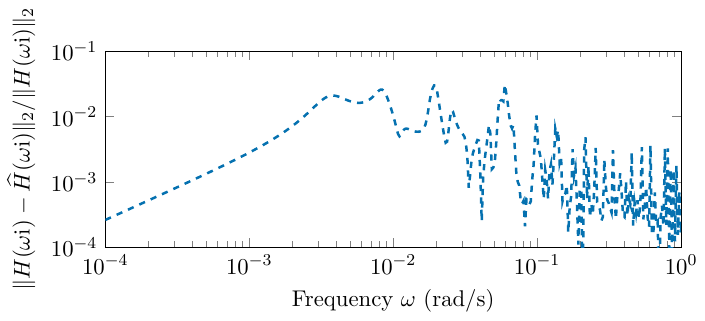}%
    {\tikzset{external/remake next}}{}%
  \input{graphics/soprbt_errorrel.tikz}%
  \tikzexternaldisable%

    \caption{Relative error.}
    \label{PBTRMV_fig:soprbt_rel}
  \end{subfigure}
  \vspace{0\baselineskip}

  \tikzexternalenable%
  \tikzsetnextfilename{soprbt_legend}%
  \filemodCmp{graphics/soprbt_legend.tikz}{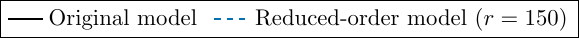}%
    {\tikzset{external/remake next}}{}%
  \begin{tikzpicture}
  \begin{axis}[
    hide axis,
    height = .75\textwidth,
    width  = .12\textheight,
    xmin   = 0,
    xmax   = 10,
    ymin   = 0,
    ymax   = 10,
    legend columns = -1,
    legend style   = {
      anchor = south, 
      /tikz/every even column/.append style={column sep = 2mm}},
        legend cell align = left
  ]
      \addlegendimage{fom}; \addlegendentry{Original model};
      \addlegendimage{rom}; \addlegendentry{Reduced-order model ($r = 150$)};
    \end{axis}
\end{tikzpicture}%
  \tikzexternaldisable%

  \caption{Positive real balanced truncation results for the triple
    chain oscillator.}
  \label{PBTRMV_fig:soprbt}
\end{figure}

Following the previously presented theory, we first compute a reduced-order
model in first-order form of order $2r = 200$ and then recover the structure
of a second-order model of order $r = 100$.
The latter has again the form
\begin{align*}
  \begin{aligned}
    \sMr \ddot{\sxr}(t) + \sDr \dot{\sxr}(t) + \sKr \sxr(t) & = \sBr \u(t), &
     \yr(t) & = \sBr^{\trans} \sxr(t),
  \end{aligned}
\end{align*}
with symmetric $\sMr, \sDr, \sKr \in \R^{r \times r}$ and
$\sBr \in \R^{r \times m}$, where
$\sMr = I_{r}$, $\sKr > 0$ and $\sDr = \sDr^{\trans}$.
The plot of the absolute value of the original and reduced-order transfer
functions can be found in \Cref{PBTRMV_fig:soprbt_tf}, whereas
\Cref{PBTRMV_fig:soprbt_rel} displays the relative error of the transfer
function on the imaginary axis, respectively.
With a maximum relative error of approximately $3.1 \cdot 10^{-2}$ we obtain a
good match between the original and the reduced-order system.

\section{\texorpdfstring{$\Hinf$}{H-infinity}-optimal rational approximation}%
\label{PBTRMV_sec:hinf_interp}

In this section, we briefly describe an interpolatory $\Hinf$ model order
reduction scheme for systems with symmetric mass, damping, and stiffness
matrices and co-located inputs and position outputs, i.e., we have
$\sM, \sD, \sK > 0$, $\sB = \sCp^{\trans}$ and $\sCv = 0$ in order to be able
to preserve symmetry and asymptotic stability by an appropriate choice of the
projection spaces.
More precisely, we construct a sequence of reduced-order transfer functions of
the form
\begin{align*}
  \begin{aligned}
    \TFr_{j}(s) & := \sB^{\trans} V_{j} \big( s^{2} V_{j}^{\trans}
      \sM V_{j} + s V_{j}^{\trans} \sD V_{j} + V_{j}^\mathsf{\trans} \sK V_{j}
      \big)^{-1} V_{j}^{\trans} \sB, & j & = 1, 2, \ldots
  \end{aligned}
\end{align*}
for appropriately chosen projection matrices $V_{j}$.
Our method aims to iteratively reduce the $\Hinf$-norm of the error transfer
function $\fE_{j}(s) := \TF(s) - \TFr_{j}(s)$.
To do so, we compute $\left\lVert \fE_j \right\rVert_{\Hinf} :=
\max_{\omega \in \R \cup \{ \infty\}} \left\lVert \fE_{j}(\omega \ri)
\right\rVert_{2}$.
The point where $\left\lVert \fE_j \right\rVert_{\Hinf}$ is attained is denoted
by $\omega_{j} \ri$.
We then choose $V_{j + 1} := \begin{bmatrix} V_{j} & (-\omega_{j}^{2} \sM +
\omega_{j} \ri \sD + \sK)^{-1} \sB \end{bmatrix}$. For numerical reasons, $V_{j+1}$ (and all other projection matrices appearing in this section) are orthogonalized.
This choice guarantees Hermite interpolation properties between the original and
reduced-order transfer functions $\TF(s)$ and $\TFr_{j}(s)$ at the interpolation
points $\omega_{1} \ri$, $\omega_{2} \ri$, $\ldots$, $\omega_{j} \ri$, such
that the error near these points becomes small.
This procedure is repeated until a specified error tolerance is met.

The main computational cost of this algorithm is the repeated computation of the
$\Hinf$-norm of the error transfer function, which is expensive to evaluate
since the error system is of large dimension.
However, these computations have been made possible by the methods presented
in~\cite{PBTRMV_AliBMetal17, PBTRMV_SchV18}, which we use here and which work
as follows.

Assume that a transfer function is given by $\TF(s) = C(sE-A)^{-1}B$ with
the regular matrix pencil $sE - A \in \R[s]^{n \times n}$, $B \in
\R^{n \times m}$, and $C \in \R^{p \times n}$ with $n \gg m, p$.
Then, the algorithm determines a sequence of reduced-order transfer functions of
the form $\TF_{j}(s) = C V_{j} (s W_{j}^{\herm} E V_{j} - W_{j}^{\herm} A
V_{j})^{-1} W_{j}^\herm B$, where $V_{j}, W_{j} \in \C^{n \times k_{j}}$ and
$k_{j} \ll n$ for $j = 1, 2, \ldots,$ and where $ \left\lVert \TF_{j}
\right\rVert_{\Linf} := \max_{\omega \in \R \cup \{\infty\}} \left\lVert
\TF_{j}(\omega \ri) \right\rVert_{2}$ converges to $\left\lVert \TF
\right\rVert_{\Hinf}$.
Since the matrices defining the reduced-order transfer functions $\TF_j(s)$,
are of small dimensions, $\left\lVert \TF_{j} \right\rVert_{\Linf}$ can be
efficiently computed using well-established methods such
as~\cite{PBTRMV_BoyB90, PBTRMV_BenSV12}.
Assume first that $m = p$.
Further suppose that $j$ interpolation points $\omega_{1} \ri, \ldots, 
\omega_{j} \ri \in \ri \R$ are already given.
In this case, we choose
\begin{align*}
  V_{j} & = \begin{bmatrix} (\omega_{1} \ri E - A)^{-1} B & \ldots &
    (\omega_{j} \ri E - A)^{-1} B \end{bmatrix}, \\
  W_{j} & = \begin{bmatrix} (\omega_{1} \ri E - A)^{-\herm} C^{\herm} & \ldots &
    (\omega_{j} \ri E - A)^{-\herm} C^{\herm} \end{bmatrix},
\end{align*}
which amounts to the Hermite interpolation conditions
\begin{align*}
  \begin{aligned}
    \TF(\omega_{k} \ri) & = \TF_{j}(\omega_{k} \ri), &
      \TF'(\omega_{k} \ri) & = \TF_{j}'(\omega_{k} \ri), &
      k & = 1, \ldots, j,
  \end{aligned}
\end{align*}
that carry over directly to the functions $\sigma(s) := \left\lVert \TF(s)
\right\rVert_{2}$ and $\sigma_{j}(s) := \left\lVert \TF_{j}(s)
\right \rVert_{2}$.
These Hermite interpolation conditions are then used to prove a superlinear rate
of convergence to a local maximum of $\sigma(\ri \cdot)$ provided that the
algorithm converges.
The situation is more difficult if $m \neq p$ since then, $V_{j}$ and $W_{j}$
would have different dimensions and the pencil $s W_{j}^{\herm} E V_{j} -
W_{j}^{\herm} A V_{j}$ would be singular.
This situation also occurs if $V_{j}$ or $W_{j}$ do not have full column rank.
Thus, an alternative choice for $V_{j}$ and $W_{j}$, which is outlined
in~\cite{PBTRMV_AliBMetal17}, and QR factorizations can be used to obtain the
matrices $V_j, W_j$ such that the pencil $s W_{j}^{\herm} E V_{j} -
W_{j}^{\herm} A V_{j}$ is regular.

\begin{figure}[t]
  \centering
  \tikzexternalenable%
  \tikzsetnextfilename{hinf_interpol}%
  \filemodCmp{graphics/hinf_interpol.tikz}{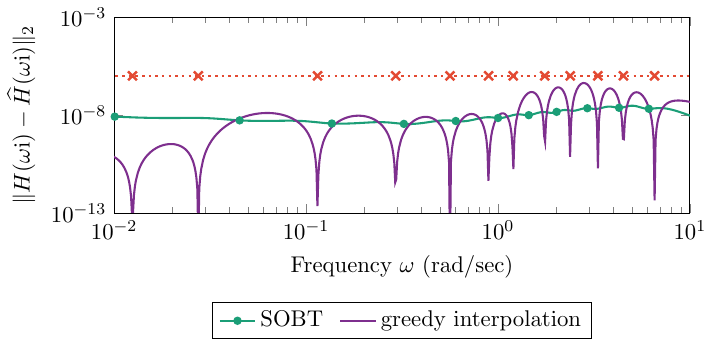}%
    {\tikzset{external/remake next}}{}%
  \input{graphics/hinf_interpol.tikz}%
  \tikzexternaldisable%

  \caption{Comparison of absolute errors of different methods for second-order
    model reduction. The interpolation points of our greedy approach are plotted
    as red crosses.}
  \label{PBTRMV_fig:hinf}
\end{figure}

In \Cref{PBTRMV_fig:hinf}, we illustrate the effectiveness of our
algorithm using the triple chain oscillator benchmark
example~\cite{PBTRMV_TruV09} (see \Cref{PBTRMV_fig:triplechain1}) of
order $n = 1\,000$ and an $\Hinf$-error bound of $10^{-6}$.
Compared with the second-order balanced truncation (SOBT, vp version) approach
described in~\cite[Alg.~3.2]{PBTRMV_ReiS08}, the greedy interpolation approach
has a slightly larger maximal error for the same reduced order $r = 29$.
However, in contrast to SOBT, we obtain full information on the current error
and may terminate whenever the reduced-order model satisfies a prescribed error
bound.
Moreover, the method also allows for an easy adaptation to frequency-limited
reduction.
We currently investigate post-processing strategies to improve the performance
of the greedy interpolation such as an additional optimization of the
interpolation points.
Furthermore, our approach may be combined with the subspaces obtained from
(SO)BT to initialize the first projection space.

The previously described algorithm leads to an error function $\fE_{j}(s)$,
which is zero at the interpolation points in exact arithmetics.
This results in the spiky shape of the error maximum singular value function
that can be observed in \Cref{PBTRMV_fig:hinf}.
Such a behavior is generally unwanted, since this indicates that our
reduced-order model approximates the given model at a few frequencies much
better than at others.
In this way, accuracy is ``wasted'' in a few regions that could be used to
improve the overall accuracy of the reduced-order model.
This is an inherent problem of a greedy approximation strategy with
interpolation on the imaginary axis.

To smoothen the error maximum singular value function and reach a better
approximation with respect to the $\Hinf$-norm, we use direct numerical
optimization.
In particular, after a new interpolation point has been chosen according to the
previously described greedy algorithm, we vary the interpolation points such
that the $\Hinf$-error is locally minimized.
This requires the solution of a nonsmooth, nonconvex and nonlinear optimization
problem.
We use the method described in~\cite{PBTRMV_CurMO17}, implemented in the
software package GRANSO\footnote{available at
\url{http://www.timmitchell.com/software/GRANSO}}.
This iterative optimization requires the repeated evaluation of the $\Hinf$-norm
of the error system, which is high-dimensional.
For that, we again apply the method described in~\cite{PBTRMV_AliBMetal17},
which is well-suited for this task.

It is important to note that the gradient of ${\lVert \fE_{j} \rVert}_{\Hinf}$
with respect to the interpolation points can be computed analytically.
Furthermore, the direct optimization strategy is not limited to just the
interpolation points.
On top of that, we can optimize \emph{tangential directions} of the
interpolation as well.
In case of tangential interpolation, the projection matrix can be chosen as
\begin{align*}
  V_{j} & = \begin{bmatrix} (s_{1}^{2} \sM + s_{1} \sD + \sK)^{-1} \sB b_{1} &
    \ldots & (s_{j}^{2} \sM + s_{j} \sD + \sK)^{-1} \sB b_{j} \end{bmatrix},
\end{align*}
where $b_{k} \in \C^{m}$ and $s_{k} \in \C$ for $k = 1,\dots, j$ are not in
the spectrum of $s^{2} \sM + s \sD + \sK$.
In this way, the interpolation condition is relaxed such that we now only have
\emph{tangential} Hermite interpolation between $\TF(s)$ and $\TFr_{j}(s)$,
that is
\begin{align*}
  \begin{aligned}
    \TF(s_{k}) b_{k} & = \TFr_{j}(s_{k})b_{k}, &
      b_{k}^{\herm} \TF(s_{k})  & = b_{k}^{\herm} \TFr_{j}(s_{k}), &
      b_{k}^{\herm} \TF'(s_{k}) b_{k} & = b_{k}^{\herm} \TFr_{j}'(s_{k})b_{k},
  \end{aligned}
\end{align*}
for $k = 1, \dots, j$.
This results in an optimization that can exploit more degrees of freedom, while
the size of the projection matrix and hence the size of the reduced-order model
is further reduced.

\begin{figure}[t]
  \centering
  \tikzexternalenable%
  \tikzsetnextfilename{hinf_errorplot2}%
  \filemodCmp{graphics/hinf_errorplot2.tikz}{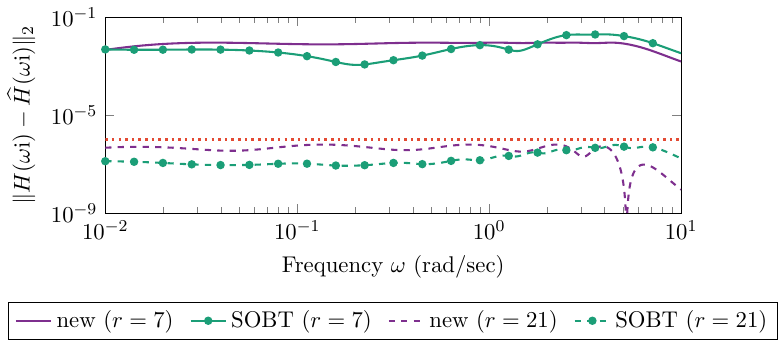}%
    {\tikzset{external/remake next}}{}%
  \input{graphics/hinf_errorplot2.tikz}%
  \tikzexternaldisable%

  \caption{Comparison of the (absolute) errors of SOBT (vp) with the new,
    optimization-based method.}
  \label{PBTRMV_fig:opt}
\end{figure}

In \Cref{PBTRMV_fig:opt}, we show a comparison of the errors between this
new method and (the faster) SOBT on the triple chain oscillator example with
$n = 301$.
The new method leads to an error function that is more steady and even
outperforms the SOBT method for the smaller reduced-order model.
However, for the slightly larger model order that is required to meet the given
error bound of $10^{-6}$, the optimization got stuck in a local optimum and
the error is less steady.
Nevertheless, the accuracy is still comparable with the accuracy obtained by
SOBT.

In order to improve the behavior of the error and make it globally more steady
we have also recently developed a new optimization framework.
There, we do not optimize the $\Hinf$-error itself but rather parametrize the
reduced-order model and then minimize the sum of squares of the error evaluated
at certain sampling points on the imaginary axis~\cite{PBTRMV_SchV20}.


\section{Conclusions}%
\label{PBTRMV_sec:conclusions}

We have presented an overview of recently developed structure-preserving
model order reduction methods for second-order systems.
We have started with an adaptation of the dominant pole algorithm for
modally damped mechanical systems and, afterwards, have introduced extensions
of the frequency- and time-limited balanced truncation methods for second-order
systems in various ways.
We have presented an approach for structure recovery of second-order systems
based on positive real balanced truncation, which also yields an a priori error
bound in the gap metric, and concluded with an $\Hinf$ greedy interpolation
approach yielding an $\Hinf$-error optimal approximation.
Numerical examples for all the presented approaches have illustrated their
effectiveness.

Additionally to the approaches for linear systems summarized here, we were able
to develop model order reduction techniques for (parametric) mechanical systems
with special nonlinearities, namely bilinear control systems and
quadratic-bilinear systems.
These techniques are described
in~\cite{PBTRMV_Wer21, PBTRMV_BenGW21, PBTRMV_BenGW21a}.


\section*{Acknowledgments}%
\addcontentsline{toc}{section}{Acknowledgments}

This work has been supported by the German Research Foundation (DFG) Priority
Program 1897: \textquotedblleft{}Calm, Smooth and Smart -- Novel Approaches
for Influencing Vibrations by Means of Deliberately Introduced
Dissipation\textquotedblright{}.
Benner and Werner were also supported by the German Research Foundation
(DFG) Research Training Group 2297
\textquotedblleft{}Mahtematical Complexity Reduction
(MathCoRe)\textquotedblright{}, Magdeburg.
  
This work has been carried out while R.~S. Beddig and M. Voigt were affiliated
with Technische Universit\"at Berlin, I. Dorschky, T. Reis, and M. Voigt were
at Universit{\"a}t Hamburg and S.~W.~R. Werner was at the Max Planck Institute
for Dynamics of Complex Technical Systems in Magdeburg. The support of these
institutions is gratefully acknowledged.


\addcontentsline{toc}{section}{References}
\bibliographystyle{plainurl}
\bibliography{bibtex/myref}

\end{document}